\documentclass[11pt,final,3p]{elsarticle}
\usepackage{amssymb}
\usepackage{amsthm}
\usepackage{amsmath}
\usepackage{latexsym}
\usepackage{amsmath}
\usepackage{xcolor}
\usepackage{graphicx}
\usepackage{indentfirst}
\usepackage[OT2,OT1]{fontenc}
\usepackage{mathrsfs}
\biboptions{numbers,sort&compress}
\usepackage{pdfsync}
\usepackage{hyperref}
\hypersetup{hypertex=true,
	colorlinks=true,
	linkcolor=blue,
	anchorcolor=blue,
	citecolor=blue}
\allowdisplaybreaks
\usepackage[utf8]{inputenc}

\renewcommand{\thefootnote}{\arabic{footnote}}

\numberwithin{equation}{section}

\newcommand\blfootnote[1]{%
	\begingroup
	\renewcommand\thefootnote{}\footnote{#1}%
	\addtocounter{footnote}{-1}%
	\endgroup
}
\journal{.}
\begin{document}
\begin{frontmatter}
\title{ The weak averaging principle of  stochastic functional partial differential equations  with H$\ddot{\text{o}}$lder continuous coefficients and infinite delay  }

\author{{ \blfootnote{$^{*}$Corresponding author } Shuaishuai Lu$^{a}$ \footnote{ E-mail address : luss23@mails.jlu.edu.cn}
		,~ Xue Yang$^{a}$}  \footnote{E-mail address : xueyang@jlu.edu.cn},~ Yong Li$^{a,b,*}$  \footnote{E-mail address : liyong@jlu.edu.cn}\\
	{$^{a}$College of Mathematics, Jilin University,} {Changchun 130012, P. R. China.}\\
	{$^{b}$Center for Mathematics and Interdisciplinary Sciences, Northeast Normal University,}
	{Changchun 130024, P. R. China.}
}

\begin{abstract}
In this paper, we establish the weak averaging principle for stochastic functional partial differential equations (in short, SFPDEs)  with H$\ddot{\text{o}}$lder  continuous  coefficients and infinite delay by a new generalized coupling approach.   Firstly, we rigorously establish the existence and uniqueness of weak solutions for a specific class of finite-dimensional systems by the generalized coupling approach. Then we extend these results to their infinite-dimensional counterparts using the variational approach and Galerkin projection technique. Subsequently, we establish the  averaging principle  for SFPDEs with infinite delay  in the weak sense, i.e.,  we prove that the solution of the original system converges in law to that of the averaged system on a finite interval $[0,T]$  as the  small parameter $\varepsilon\to 0$. To illustrate our findings, we present two applications: stochastic generalized porous media equations and stochastic reaction-diffusion equations.
~\\
~\\
 \textbf{keywords}: Weak averaging principle; H$\ddot{\text{o}}$lder continuous  coefficients; Stochastic functional partial differential equations; Infinite delay; Generalized coupling
 ~
\end{abstract}
\end{frontmatter}
\section{\textup{Introduction}}
	Random phenomena encompass both natural occurrences and man-made systems, ranging from  financial market fluctuations to molecular motions in living organisms, all of which can be mathematically described by stochastic differential equations. Traditional differential equation theory is rooted in deterministic assumptions. However, many real-world systems are influenced by external random disturbances. Therefore, the introduction of  stochastic processes and stochastic differential equations becomes essential for accurately modeling and analyzing these systems. Practical applications of stochastic models, such as climate interactions \cite{ref7,ref8}, financial fluctuations \cite{ref9,ref10} and geophysical fluid dynamics \cite{ref11,ref12}, frequently  exhibit significant oscillatory components. These high-frequency oscillations present challenges for direct analysis and simulation of the system's properties. Consequently, deriving simplified equations that effectively capture the long-term evolution of these systems is crucial. The key strategy involves applying the averaging principle, which  ``averages out" highly oscillatory components under suitable conditions, thereby yielding an averaged system that is more manageable for analysis.   This averaged system governs the evolution of the original complex system over extended time periods, enabling us to grasp its fundamental dynamics independently of high-frequency oscillations.

Consider  the following stochastic functional partial differential equations (SFPDEs)  with high-frequency oscillating variables and infinite delay on a separable Hilbert space $U_{1}$:
\begin{align}\label{j1}
 \text{d}u^{\varepsilon }(t)=[A(u^{\varepsilon }(t))+f(\frac{t}{\varepsilon } ,u^{\varepsilon }_{t})]\text{d}t+g(\frac{t}{\varepsilon } ,u^{\varepsilon }_{t})\text{d}W(t),
	\end{align}
and
\begin{align}\label{j2}
 \text{d}u(t)=[A(u(t))+f^{*}(u_{t})]\text{d}t+g^{*}(u_{t})\text{d}W(t),
	\end{align}
where $\varepsilon\in (0,1)$. The function $f^{*}$, $g^{*}$ are defined as  $$f^{*}(\varphi )=\lim_{T \to \infty} \frac{1}{T}\int_{0}^{T} f(t,\varphi ) \text{d}t,\quad \lim_{T \to \infty}\frac{1}{T}  \int_{0}^{T} \left \| g(t,\varphi )-g^{*}(\varphi ) \right \|^{2}_{\mathscr{L}(U_{2},U_{1})} \text{d}t=0$$
uniformly with respect to $t\in \mathbb{R}^{+} $ for any $\varphi \in C((-\infty , 0],U_{1})$. Let $W(t)$ denote a  cylindrical Wiener process  on a separable Hilbert space $U_{2}$ and let $u_{t} =\{u_{t} (\theta )\}=u(t+\theta )(-\infty  \le \theta \le 0)$ be the segment process or solution map of system. A fundamental question arises: Does the solution of system \eqref{j1} converge to that of the averaged system \eqref{j2} as the  small parameter $\varepsilon\to 0 $? This question arises naturally from both physical and mathematical perspectives, and the theory addressing such problems is the averaging principle mentioned above.

The concept of the averaging principle originated from the study of   nonlinear oscillations of deterministic systems by Bogoliubov \cite{ref13,ref16}, and its extension to stochastic differential equations was further advanced  by Khasminskii \cite{ref15}. The averaging principle for finite dimensional systems has been extensively studied in recent decades, building on the pioneering work of Khasminskii; see, for example, \cite{ref19,ref20,ref21,ref22,ref23,ref25,ref26,ref27,ref28,ref38}. For infinite-dimensional systems, the averaging principle was introduced by Henry \cite{ref17}. Subsequently, the corresponding averaging principle for infinite-dimensional systems has also undergone  significant development. In \cite{ref18}, Maslowski et al. presented an averaging principle for stochastic evolution equations with small parameters using the semigroup method. Kuksin and Piatnitski \cite{ref29} explored the Whitham averaging method in damped-driven equations affected by random disturbances. Cerrai and Freidlin \cite{ref31} demonstrated the averaging principle for stochastic reaction-diffusion equations where the diffusion coefficients and the rates of reactions have different orders.  The averaging principle for nonautonomous stochastic reaction-diffusion equations exhibiting slow-fast dynamics was addressed in \cite{ref32}. Furthermore, Cheng and Liu \cite{ref2} established the second Bogolyubov theorem and a global averaging principle for a specific class of SPDEs with monotone coefficients. For further details on this topic, we refer the interested reader to \cite{ref36,ref37,ref39,ref50,ref51,ref52,ref30,ref33}.

The averaging principle fundamentally requires Lipschitz regularity of coefficients---a condition ensuring solution uniqueness and continuous dependence on initial data. Additionally, it allows us to work within the framework of square-integrable functions, thereby simplifying the proof process for a priori bounds of solutions in both fast and slow systems. Despite its importance in studying the dynamic properties of stochastic systems, the Lipschitz continuity condition is often not satisfied in many significant stochastic models. For instance,  the Cox-Ingersoll-Ross model and diffusion coefficients in the Feller branch diffusion exhibit only H$\ddot{\text{o}}$lder continuous, not Lipschitz continuous. This limitation restricts the applicability of the averaging principle in various relevant scenarios.  Consequently, there has been increasing interest in recent years in developing averaging principles of stochastic models with coefficients that are not Lipschitz continuous. For example, Veretennikov \cite{ref28} investigated the averaging principle or SDEs with slow and fast time scales. Specifically, the drift coefficients of the equations involving the slow variables must be bounded and measurable, while all other coefficients are assumed to be globally Lipschitz continuous. In \cite{ref55}, Cerrai presented  the averaging principle for a system of slow and fast reaction-diffusion equations featuring multiplicative noise and coefficients with polynomial growth. Sun et al. \cite{ref56} established the averaging principle for slow-fast SPDEs with H$\ddot{\text{o}}$lder continuous drift coefficients. In \cite{ref57}, R$\ddot{\text{o}}$ckner et al.  studied the averaging principle for semi-linear slow-fast SPDEs featuring additive noise with reaction coefficients that are merely assumed to be H$\ddot{\text{o}}$lder continuous with respect to the fast variable. For systems similar to system \eqref{j1} that involve high-frequency oscillating variables, Xu and Xu \cite{ref58} investigated the averaging principle of stochastic evolution equations under non-Lipschitz conditions, within  the framework of mild solutions. In addition,  for the  diffusion approximations of SDEs with singular coefficients, we refer interested reader to \cite{ref59,ref60,ref61}.

The current issue pertains to whether the aforementioned averaging principle can be established for a stochastic differential system with parameters exhibiting high-frequency oscillations, assuming the diffusion and drift coefficients are H$\ddot{\text{o}}$lder  continuous. There are limited studies on this question, and our article addresses the critical problem. We will explore the averaging principle under H$\ddot{\text{o}}$lder conditions within the framework of stochastic partial differential equations with infinite delay. However, we emphasize that the conclusions in the paper are also applicable to corresponding models with finite delay and without delay, with slight adjustments to the conditions. In fact, the derivation process often becomes more straightforward in these cases. The reason we have chosen to focus on system \eqref{j1} for our study is twofold: it allows us to showcase our methodology effectively and is also widely applicable in various contexts. Since time delays are pervasive  in many real-world systems, they play a crucial role in mathematical and physical models. The motivation arises from phenomena involving delayed transmission, such as high-speed fields in wind funnel experiments, species growth patterns, and the incubation periods in disease models. In these systems, time delays and high-frequency oscillations play crucial roles in determining the behavior of the system over time. For example, in wind funnel experiments, the high-frequency fluctuations in air speed can lead to complex dynamic behaviors that necessitate the use of stochastic models to capture the full range of possible outcomes. Similarly, the delayed transmission in disease models can significantly affect the spread of the disease, making the incorporation of delay and oscillations in the model essential for accurate predictions. The theory of stochastic functional differential equations with infinite delays has been garnering increasing attention, as reflected in recent literature; for examples, see \cite{ref1,ref62,ref63,ref65,ref66} and others.

Establishing solution well-posedness constitutes the foundational requirement for developing the averaging principle. Kulik and Scheutzow \cite{ref68} pioneered weak uniqueness through generalized coupling for finite-dimensional SFDEs  with H$\ddot{\text{o}}$lder continuous coefficients. Subsequently, Han further \cite{ref777} established the well-posedness and small mass limit for the stochastic wave equation with H$\ddot{\text{o}}$lder continuous noise coefficients. The well-posedness theory for non-Lipschitz systems has been extensively developed in \cite{ref62,ref67,ref69,ref70}. In this paper, we  analyze the well-posedness of solutions to the system,  inspired by the generalized coupling method proposed in \cite{ref68}. For SFDEs with H$\ddot{\text{o}}$lder  continuous coefficients, we first prove  the existence and uniqueness of weak solutions for a class of finite-dimensional systems(see Theorem 3.1). Subsequently, in the framework of infinite-dimensional systems, we employ the Galerkin projection technique. By combining this method with insights from finite-dimensional analysis, we establish the existence and uniqueness of probabilistic weak solutions, as stated in Theorem 3.3.  The primary challenge lies in proving the uniqueness of weak solutions.  Under Lipschitz conditions,  the dependence of solutions on initial data can be analyzed via Gronwall's lemma; however, for  H$\ddot{\text{o}}$lder continuous coefficients, this approach is not realistic. To address this limitation, we develop a novel probabilistic framework based on generalized coupling, which ensures uniqueness  in law for weak solutions(i.e., the weak solution is unique in law). In this paper, we provide a comprehensive overview of the study of S(F)PDEs with H$\ddot{\text{o}}$lder continuous diffusion and drift coefficients in an abstract framework, with applications extendable to stochastic models with low-regularity coefficients.

Subsequently, another major result in the present paper is  to establish the weak averaging principle for systems \eqref{j1} with H$\ddot{\text{o}}$lder  continuity coefficients. Namely, we prove that the convergence of solutions from the original Cauchy problem to those of the averaged equation across the finite interval $[0, T]$. In fact, under Lipschitz continuity conditions, convergence in the mean-square sense can be obtained, i.e.:
\begin{align*}
		\lim_{\varepsilon  \to 0} \mathbb{E}\underset{t\in [0,T]}{\sup}   \left \|u^{\varepsilon }(t;\varphi ^{\varepsilon }) -u^{* }(t;\varphi ^{*}) \right \|^{2}_{U_{1}}=0.
	\end{align*}
However, when the coefficients satisfy the H$\ddot{\text{o}}$lder continuity condition, mean-square convergence does not hold, and weaker probabilistic topologies are necessitated. To this end, convergence in law is established via the generalized coupling method and the classical Khasminskii time discretization. That is why we call it ``weak averaging principle''. Our main convergence result states: For any $T > 0$, with initial data satisfying $\lim_{\varepsilon  \to 0} \mathbb{E} \left \|\varphi ^{\varepsilon } -\varphi ^{*} \right \|^{2}_{h}=0$, we establish
\begin{align}\label{c1}
		\lim_{\varepsilon  \to 0}|\mathbb{E}\mathcal{L}(u^{\varepsilon }(t;\varphi ^{\varepsilon })| _{[0,T]}) - \mathbb{E}\mathcal{L}(u^{* }(t;\varphi ^{*})|_{[0,T]})|= 0,
	\end{align}
 where $\mathcal{L}$ is an arbitrary real-valued, bounded  continuous function defined on $C([0, T], U_{1})$, $u^{\varepsilon }(s;\varphi ^{\varepsilon })$ is the weak solution of \eqref{j1} with the initial value $u^{\varepsilon }_{0}=\varphi ^{\varepsilon }$ and
$u^{*}(s;\varphi ^{*})$ is the weak solution of \eqref{j2} with the initial value $u^{*}_{0}=\varphi ^{*}$. To the best of our knowledge, there is no work so far on averaging principle  for  S(F)PDEs with H$\ddot{\text{o}}$lder continuous diffusion and drift coefficients. Unlike classical results relying on Lipschitz conditions, our framework enables the application of the averaging principle to stochastic models with low-regularity coefficients.

The remainder of this paper is organized as follows. In Section 2, we introduce the assumptions and state our main results. Section 3 rigorously establishes the existence and uniqueness of weak solutions for a class of finite-dimensional systems with H$\ddot{\text{o}}$lder  continuous drift and diffusion coefficients. Subsequently, employing Galerkin projection techniques and insights from finite-dimensional systems, we extend these results to demonstrate the existence and uniqueness of probabilistic weak  solutions for corresponding infinite-dimensional systems. In Section 4, we concentrate on proving the averaging principle under the condition that the drift and diffusion coefficients of system \eqref{j1} satisfy H$\ddot{\text{o}}$lder  continuity. Additionally, we also  establish the averaging principle for systems with both finite delay and delay-free under the H$\ddot{\text{o}}$lder  continuity condition. In Section 5, we will give two examples to illustrate the applicability of our result.

	\section{\textup{Preliminaries}}
	The $(U_{i},\left \| \cdot  \right \| _{U_{i}})$, $i=1,2$ are  separable Hilbert spaces with  inner product $ \left \langle \cdot, \cdot \right \rangle _{U_{i}}$. The space
$(B,\left \| \cdot  \right \| _{B})$ is a reflexive Banach space and  $U_{i}^{*}$, $B^{*}$ denote the dual spaces of  $U_{i}$, $B$, respectively. Let
\begin{align*}
		B\subset  U_{1}\subset B^{*},
	\end{align*}
where the embedding $B\subset U_{1}$ is continuous and dense.  Thus,
$U_{1}^{*}$ is densely and continuously embedded in $B ^{*}$. Denote the pairing between $B^{\ast}$ and $B$ by $_{B ^{\ast}} \langle \cdot , \cdot \rangle _{B} $, which implies that for all $u \in U_{1}$, $v \in B$,
\begin{align*}
		_{B ^{\ast}} \langle u, v\rangle _{B} = \langle u, v\rangle_{ U_{1}},
	\end{align*}
 and  ($B,U_{1}, B ^{\ast}$) is called a Gelfand triple. The  $\left ( \Omega ,\mathscr{F},\mathbb{P} \right )$ is a certain complete probability space with a filtration $\lbrace\mathscr{F}_{t}\rbrace_{t\ge0}$ satisfying the usual condition and   $\mathbb{R}^{n}$ is the $n$-dimensional Euclidean space equipped with the norm $\left | \cdot  \right | $. If $K$ is a matrix or a vector, $K^{'}$ is its transpose. For a matrix $K$, the norm is expressed as $\left \| K\right \|= \sqrt{\text{trace}(KK^{'})}$ and denote by $\left \langle \cdot , \cdot  \right \rangle $ the inner product of $ \mathbb{R}^{n} $.  For any $q \ge 1$  and  the Banach space $(X, \left \|   \cdot  \right \|_{X}   )$,   $C_{b}(\mathbb{R}, L^{q}(\mathbb{P},X)) $ to represent the set of all continuous and uniformly bounded stochastic processes from $\mathbb{R}$ into $L^{q}(\mathbb{P},X)$. If $(X,\left \| \cdot \right \| _{X}  )$ is  a real separable Hilbert spaces, for a given $h > 0$, we give the following definition of spaces $C^{h}$:
     \begin{eqnarray*}
		C^{h}_{X} =\left \{ \varphi \in C((-\infty ,0 ],X ):\underset{\theta \to -\infty }{\lim }e^{h\theta }\varphi (\theta )=x  \in   X \right \}.
	\end{eqnarray*}
 The space $C^{h}_{X}$  is a Banach space with norm $\left \| \varphi  \right \| _{h } =\underset{-\infty < \theta \le 0}{\sup}\left \| e^{h\theta }\varphi (\theta ) \right \| _{X}  $, and $C^{h}_{X}$ has the following properties(\cite{ref1}):
\begin{enumerate}[(1)]
		\item For any $T>0$, $x(\cdot):(-\infty ,T]\to X$ is continuous on $[0, T )$ and $x_{0}=\{x(\cdot)\}:(-\infty ,0 ]\to  X $ is in $C^{h}_{X}$. Then for every $t\in  [0, T ]$,
        \begin{enumerate}[(a)]
		\item $x_{t}(\cdot ) =\{x(t+\cdot )\}\in C^{h}_{X}$;
	    \end{enumerate}
\begin{enumerate}[(b)]
		\item $\left \| x(t) \right \| _{X} \le \left \| x_{t}  \right \| _{h } $.
	    \end{enumerate}
\item  $x_{t}$ is a $C^{h}_{X}$-valued continuous function for any $t\in [0, T ]$.
		\item The space $C^{h}_{X}$ is a complete  space.
	\end{enumerate}

Let $\mathcal{B}(X)$ denote the $\sigma $-algebra generated by space $X$
and $\mathcal{P}(X)$ be the family of all probability measures on $(X , \mathcal{B}(X))$ with the following metric
 \begin{align*}
		d_{X}( \mu ,\nu ) :=\sup\{\left | \int \Upsilon \text{d}\mu _{1}- \int \Upsilon \text{d}\mu _{2} \right |:\left \| \Upsilon  \right \|_{BL}\le 1  \},
	\end{align*}
 where $\left \| \Upsilon  \right \|_{BL}:=\left \| \Upsilon  \right \| _{\infty }+Lip(\Upsilon )$ and $\left \|\Upsilon  \right \| _{\infty }=\sup_{\varphi \in X}\frac{|\Upsilon (\varphi)|}{(1+\left \| \varphi  \right \|_{X} )^{2}} $, $Lip(\Upsilon )=\sup_{\varphi_{1}\ne \varphi_{2}}\frac{\left | \Upsilon (\varphi_{1})-\Upsilon (\varphi_{2}) \right | }{\left \|\varphi_{1}-\varphi_{2}  \right \|_{X}} $. Then it is not difficult to verify that the space  $(\mathcal{P}(X), d_{X})$ is a complete metric space.

Consider the following SFPDEs with infinite delay:
\begin{align}\label{r1}
		\begin{cases}
 \text{d}u(t)=(A(t,u(t))+f(t,u_{t})\text{d}t+g(t,u_{t})\text{d}W(t), \\
u_{0} =\varphi\in C^{h}_{U_{1}},
\end{cases}
	\end{align}
 where $W(t)$ be a  cylindrical Wiener processes  on a separable Hilbert space $(U_{2},\langle \cdot , \cdot \rangle_{ U_{2} })$ with respect to a complete
filtered probability space $\left ( \Omega ,\mathscr{F},\mathbb{P} \right )$, and $u_{t} =\{u_{t} (\theta )\}=u(t+\theta )(-\infty  \le \theta \le 0)$ is the segment process or solution map of system \eqref{r1}. The measurable maps: $$A:\mathbb{R} ^{+}\times B\to B^{*},\quad f:\mathbb{R} ^{+}\times C^{h}_{U_{1}} \to U_{1},\quad g:\mathbb{R} ^{+}\times C^{h}_{U_{1}}  \to \mathscr{L}(U_{2},U_{1}),$$
where $\mathscr{L}(U_{2},U_{1})$ is  the space of all Hilbert-Schmidt operators from $U_{2}$ into $U_{1}$.

In the next section, the following lemma plays a crucial role.
\\\textbf{Lemma 2.1.} (Lemma B.1, \cite{ref68}) Let $V(t) \geq 0$ be an It\^o process with
\begin{align*}
    dV(t) = \eta(t) dt + dM(t),
\end{align*}
where $M$ is a continuous local martingale with quadratic variation
\begin{align*}
    \langle M \rangle (t) = \int_0^t m(s) ds, \quad t \geq 0.
\end{align*}
Assume that for some constants $A \geq 0$, $B > 0$, $\lambda > 0$ and a random variable $\varsigma \geq 0$
\begin{align*}
    \eta(t) \leq -\lambda V(t) + A, \quad m(t) \leq B, \quad t \leq \varsigma.
\end{align*}
Assume also that $\varsigma \leq T$ for some constant $T > 0$.
Then for every $\delta \in (0, 1/2)$ there exist constants $C_1, C_2 > 0$, which depend only on $\delta$ and $T$, such that
\begin{align*}
    P \left( \sup_{t \leq \varsigma} \left( V(t) - e^{-\lambda t} V(0) \right) \geq A \lambda^{-1} + B^{1/2} \lambda^{-\delta} R \right) \leq C_1 e^{-C_2 R^2}, \quad R \geq 0.
\end{align*}

\section{\textup{Existence and uniqueness }}
The investigation into the existence and uniqueness of solutions to equation \eqref{r1} with H$\ddot{\text{o}} $lder continuous coefficients, assumes that the initial value  \(\varphi \in C^{h}_{U_{1}}\) is  independent of  \(W(t)\). To begin our analysis, we assume that the coefficients in equation \eqref{r1} satisfy the following hypotheses:
 ~\\
 \textbf{(H1)} (Continuity) For all $t\in\mathbb{R} ^{+}$ and $u,v \in B $, the map
\begin{align*}
		\mathbb{R} ^{+}\times B\ni (t,u)\to_{B ^{\ast}}\langle A(t,u),v\rangle _{B}
	\end{align*}
is continuous.
~\\
\\\textbf{(H2)} (Growth) For $A$, there exist constant $\alpha_{1}, M>0$ for all $t\in\mathbb{R} ^{+}$  and $u\in B$ such that
 \begin{align*}
		\left \| A(t,u)\right \|_{B^{*}}^{\frac{p}{p-1}}   \le \alpha  _{1}\left \| u\right \|_{B}^{p}+ M,
	\end{align*}
and for the continuous functions $f$, $g$, there exist constant $\alpha _{1}$ and $M$ for all $t\in\mathbb{R} ^{+}$ and $\varphi\in C^{h}_{U_{1}} $ such that
\begin{align*}
		\left \| f(t,\varphi) \right \|_{U_{1}} \vee \left \| g(t,\varphi ) \right \| _{\mathscr{L}(U_{2},U_{1})}  \le \alpha _{1}\left \| \varphi\right \|_{h  }+M.
	\end{align*}
\\\textbf{(H3)} (Coercivity)  There exist constant  $ \alpha _{1}$, $M$, $p\ge2$ and $\alpha _{2} \in \mathbb{R}$ such that for all $t\in\mathbb{R} ^{+}$ and $u \in B $
\begin{align*}
		_{B ^{\ast}} \langle A(t,u), u\rangle _{B} \le-\alpha _{1}\left \| u \right \|_{B}^{p}+ \alpha_{2}\left \| u \right \| _{U_{1}}^{2}+M.
	\end{align*}
\\\textbf{(H4)}  There exist constants $\beta \in(0,1] $ and $\gamma  \in(\frac{1}{2},1 ]$. The map $A$ satisfies, for all $t\in\mathbb{R} ^{+}$ and $u,v\in B$
\begin{align*}
		2_{B ^{\ast}} \langle A(t,u)-A(t,v), u-v\rangle _{B} \le \alpha _{1}\left \| u-v \right \| _{U_{1}}^{2},
	\end{align*}
and   the functions $f$, $g$ satisfy,  for all $t\in\mathbb{R} ^{+}$ and $\varphi,\phi\in C^{h}_{U_{1}}  $ with $\left \| \varphi\right \| _{h}\vee\left \| \phi  \right \| _{h}\le M$,
\begin{align*}
\left \langle f(t,\varphi)-f(t,\phi),\varphi(0)-\phi(0) \right \rangle _{U_{1}}\le L_{M}\left \| \varphi -\phi  \right \|^{\beta+1  } _{h}
\end{align*}
and
\begin{align*}
		\left \| g(t,\varphi )- g(t,\phi ) \right \| _{\mathscr{L}(U_{2},U_{1})} \le L_{M}\left \| \varphi -\phi  \right \|^{\gamma  } _{h}.
	\end{align*}
\\\textbf{(H5)} For each \( u \in U_{1} \), \( g \) has a right inverse \( g^{-1} \) on \(U_{1}\) in the sense that
\[
g [g^{-1} u] = u \quad \text{for any } u \in U_{1},
\]
with \( g^{-1}[U_{1} ]\subset U_{2} \) and satisfies,
\[
\sup_{\varphi \in C^{h}_{U_{1}}} \left\| g^{-1}(t, \varphi) \right\|_{\mathscr{L}(U_{1}, U_{2})} < \infty,
\]
for any $t>0$.

In this paper, $L_{M}$ denotes certain positive constants dependent on $M$. To prevent ambiguity, we maintain the assumption throughout that the constants \(\alpha_1>0 , M>0\) and \(\alpha_2\in \mathbb{R}\),which are allowed to vary from line to line.

Extensive research has focused on the properties of existence and uniqueness of solutions when the system coefficients satisfy Lipschitz conditions. However, for coefficients with lower regularity, such as H$\ddot{\text{o}} $lder continuity, research on the topic is relatively scarce.  We investigate the existence and uniqueness of solutions, as well as the asymptotic behavior  of equation \eqref{r1} using the Galerkin-type approximation technique. To enhance clarity, we initially analyze the following finite-dimensional system with H$\ddot{\text{o}} $lder continuous coefficients:
\begin{align}\label{r3} \text{d}x(t)=F(t,x_{t})\text{d}t+G(t,x_{t})\text{d}B(t),
	\end{align}
where the initial data $x_{0} =\varphi \in C^{h}_{\mathbb{R}^{n} }$, $B(t)$ is an $m$-dimensional Wiener process and $F:\mathbb{R} ^{+}\times C^{h}_{\mathbb{R}^{n} }  \to \mathbb{R}^{n}$, $G:\mathbb{R} ^{+}\times C^{h}_{\mathbb{R}^{n} } \to \mathbb{R}^{n\times m}$ are two continuous maps. We assume that the coefficients in equation \eqref{r3} satisfy the following hypotheses:
~\\
\\\textbf{(h1)} The functions $F$ and $G$ are continuous in $(t,\varphi )$ and satisfy, for all $t\in\mathbb{R} ^{+}$ and $\varphi \in C^{h}_{\mathbb{R}^{n} } $
 \begin{eqnarray*}
		\left \langle F(t,\varphi) ,\varphi(0) \right \rangle   \vee \left \| G(t,\varphi ) \right \|^{2} \le \alpha _{1}\left \| \varphi\right \|_{h  }^{2}+M.
	\end{eqnarray*}
\textbf{(h2)}  The functions $F,G$ satisfy, for all $t\in\mathbb{R} ^{+}$ and $\varphi , \phi \in C^{h}_{\mathbb{R}^{n} } $  with $\left \| \varphi\right \|_{h}\vee\left \| \phi \right \|  _{h}\le M$,
\begin{align*}
		\left \langle F(t,\varphi)-F(t,\phi),\varphi(0)-\phi(0) \right \rangle\le \left \| \varphi -\phi  \right \|^{\beta+1  } _{h}
	\end{align*}
and
\begin{align*}
		 \left \| G(t,\varphi )- G(t,\phi ) \right \|  \le L_{M}\left \| \varphi -\phi  \right \|^{\gamma  } _{h}.
	\end{align*}
\textbf{(h3)} For each $\varphi\in C^{h}_{\mathbb{R}^{n}}$, there exists a right inverse $G^{-1}$ of the coefficient $G$, and
\[
\sup_{\varphi \in C^{h}_{\mathbb{R}^{n}}} \left\| G^{-1}(t, \varphi) \right\|< \infty,
\]
for any $t>0$.

~\\
\textbf{Theorem 3.1.} \emph{Consider \eqref{r3}. Suppose that the assumptions \textbf{(h1)}$-$\textbf{(h3)} hold. Then the following statement holds: for any $\varphi \in C^{h}_{\mathbb{R}^{n} } $,  system \eqref{r3} has  a probabilistic weak solution  $u(t)$ and a segment process $x_{t}$ to \eqref{r3} with $x_{0}=\varphi$,  and the weak solution $u(t)$  is unique in law, i.e., any two such solutions with the same initial value have the same law.
~\\
\\\textbf{proof:}} Let the support of $\rho \in C^{\infty }(\mathbb{R} )$ be contained in $\mathcal{O} ^{\mathcal{M} }=\{\varphi \in C^{h}_{\mathbb{R}^{n} }:\left \| \varphi  \right \| _{h}\le \mathcal{M} \}$, i.e.,
\begin{align*}
\left\{\begin{matrix}
 \rho (\left \| \varphi  \right \| _{h} )=1, & \varphi \in \mathcal{O} ^{\mathcal{M} }, \\
   \rho (\left \| \varphi  \right \| _{h} )=0, & \varphi \notin  \mathcal{O} ^{\mathcal{M} }.&
\end{matrix}\right.
\end{align*}
For any $ \varphi\in C^{h}_{\mathbb{R}^{n} }$, let $F^{\mathcal{M} } (t,\varphi)=F (t,\varphi)\rho (\left \| \varphi \right \|_{h} )$ and $G^{\mathcal{M} } (t,\varphi)=G (t,\varphi)\rho (\left \| \varphi \right \|_{h} )$, and we consider the following system:
\begin{align}\label{1r1} \text{d}x^{n}(t)=F^{\mathcal{M}}(t,x^{n}_{t})\text{d}t+G^{\mathcal{M}}(t,x^{n}_{t})\text{d}B(t).
	\end{align}
where $F^{\mathcal{M} } (t,\varphi)$, $G^{\mathcal{M} } (t,\varphi)$ are uniformly bounded on $ C^{h}_{\mathbb{R}^{n} }$ and  H$\ddot{\text{o}}$lder continuous.
 Thus, there exist two sequences of uniformly bounded local Lipschitz continuous functions $F^{\mathcal{M}}_{n} (t,\varphi)$ and $G^{\mathcal{M} }_{m} (t,\varphi)$ such that
 \begin{align}\label{t9}
\left | F^{\mathcal{M}}_{n} (t,\varphi)-F^{\mathcal{M}} (t,\varphi) \right | \to 0,\quad \left \| G^{\mathcal{M}}_{n} (t,\varphi)-G^{\mathcal{M}} (t,\varphi) \right \| \to 0,\quad \text{as}\quad n\to\infty,
\end{align}
 uniformly on each compact subset of $C^{h}_{\mathbb{R}^{n} }$ and $t\in[0,T]$.
Then  we further consider the following equation:
\begin{align}\label{r5} \text{d}x^{n}(t)=F^{\mathcal{M}}_{n}(t,x^{n}_{t})\text{d}t+G^{\mathcal{M}}_{n}(t,x^{n}_{t})\text{d}B(t),
	\end{align}
where $n\ge 1$. Based on the previous analysis, it is established that \eqref{r5} admits a unique strong solution \( x^{n}(t) \) with the initial condition \( x^{n}_{0} = \varphi \).

Our main strategy involves initially establishing the existence of weak solutions for \eqref{1r1} through the construction of strong solutions $x^{n}(t)$ for \eqref{r5}. Subsequently, we proceed to prove the  uniqueness in law of these weak solutions.

For any fixed \( q \ge 2 \), applying It$\hat{\text{o}} $ formula to \( |x^{n}(t)|^{q} \), we derive
\begin{align*}
\begin{split}
|x^{n}(t)|^{q}&=\left \| \varphi  \right \| _{h }^{q}+\int_{0}^{t} [\frac{q(q-1)}{2} \left | x^{n}(s) \right |^{q-2} \left \| G^{\mathcal{M}}_{n}(s,x^{n}_{s}) \right \|^{2} \\&~~~+q\left | x^{n}(s) \right |^{q-2}\left \langle F^{\mathcal{M}}_{n}(s,x^{n}_{s}),x^{n}(s)  \right \rangle ]\text{d}s
\\&~~~+\int_{0}^{t}q\left | x^{n}(s) \right |^{q-2} ( x^{n}(s))' G^{\mathcal{M}}_{n}(s,x^{n}_{s})\text{d}B(s).
\end{split}
\end{align*}
There exist constants \( L_{q}, M_{q} \) by Young's inequality and assumptions \textbf{(h1)}, such that
\begin{align*}
\begin{split}
|x^{n}(t)|^{q}&=\left \| \varphi  \right \| _{h }^{q}+L_{q}\int_{0}^{t} (\left \| x^{n}_{s} \right \| _{h}^{q}+M_{q})\text{d}s
+L_{q}\int_{0}^{t}(\left \| x^{n}_{s} \right \| _{\mathcal{D} }^{q}+M_{q})\text{d}B(s).
\end{split}
\end{align*}
  Note
\begin{align}\label{z1}
\left \| x^{n}_{s} \right \| _{h}&=\underset{-\infty < \theta \le 0}{\sup}\left | e^{h\theta } x^{n} (s+\theta ) \right | \nonumber \\&=\underset{-\infty < r \le s}{\sup}\left | e^{h(r-s) } x^{n} (r) \right |
\\&\le\underset{-\infty < r \le 0}{\sup}\left | e^{h(r-s) } x^{n} (r) \right |+\underset{0 \le r \le s}{\sup}\left | e^{h(r-s) } x^{n} (r) \right | \nonumber
\\&\le\left \| \varphi  \right \|  _{h}+\underset{r\in[0,s]}{\sup}\left | x^{n}(r) \right |. \nonumber
\end{align}
By Cauchy's inequality, Burkholder-Davis-Gundy inequality and \eqref{z1}, we have
\begin{align*}
\mathbb{E} \underset{z\in[0,t]}{\sup}  \left | x^{n}(t) \right | ^{2q}&\le L_{q} \left \| \varphi  \right \| _{h }^{2q}+L_{q}\mathbb{E}[\int_{0}^{t} (\underset{r\in[0,s]}{\sup}\left | x^{n}(r) \right |^{q}+M_{q})\text{d}s]^{2}
\\&~~~+L_{q}\mathbb{E}(\underset{z\in[0,t]}{\sup}\int_{0}^{t}(\left \| x^{n}_{s} \right \| _{\mathcal{D} }^{q}+M_{q})\text{d}B(s)) ^{2}
\\&\le L_{q} \left \| \varphi  \right \| _{h}^{2q}+L_{q}\mathbb{E}\int_{0}^{t} (\underset{r\in[0,s]}{\sup}\left | x^{n}(r) \right |^{2q}+M_{q})\text{d}s
+L_{q}\int_{0}^{t} [\mathbb{E}\left \| x^{n}_{s} \right \| _{h }^{2q}+M_{q}]\text{d}s
\\&\le L_{q}\left \| \varphi  \right \| _{h }^{2q}+L_{q}[\int_{0}^{t} \mathbb{E}\underset{z\in[0,s]}{\sup}  \left | x^{n}(z) \right | ^{2q}\text{d}s+M_{q}t].
\end{align*}
Applying the Gronwall's inequality then gives
\begin{align}\label{r6}
\begin{split}
\mathbb{E} \underset{z\in[0,t]}{\sup}  \left | x^{n}(t) \right | ^{2q}\le L_{q,M}(\left \| \varphi  \right \| _{h }^{2q}+t+e^{t})<\infty.
\end{split}
	\end{align}
For any \( t \in [0,T] \) and \( T > 0 \), we deduce from condition \textbf{(h1)} that \( F^{\mathcal{M}}_{n} \) and \( G^{\mathcal{M}}_{n} \) are bounded on every bounded subset  of \( C^{h}_{\mathbb{R}^{n} } \). Therefore, by \textbf{(h1)} and \eqref{r6}, there exists a constant \( L_{T} \), independent of \( n \), such that
\begin{align*}
\left |F^{\mathcal{M}}_{n}(s,x^{n}_{s} )\right | \vee \left \|  G^{\mathcal{M}}_{n}(s,x^{n}_{s})\right \|\le L_{T}.
\end{align*}
Hence for any $0\le z, t \le T < \infty  $, we have
\begin{align}\label{r7}
\begin{split}
&\underset{n\ge1}{\sup}\mathbb{E} \left | x^{n}(t)-x^{n}(z) \right | ^{2q}\\&\le L_{q}\underset{n\ge1}{\sup}\mathbb{E} \left | \int_{z}^{t} F^{\mathcal{M}}_{n}(s,x^{n}_{s})\text{d}s \right | ^{2q}
+L_{q}\mathbb{E} \underset{n\ge1}{\sup}\left | \int_{z}^{t} G^{\mathcal{M}}_{n}(s,x^{n}_{s})\text{d}B(s) \right | ^{2q}
\\&\le L_{q,T}\left | t-z \right | ^{q}.
\end{split}
	\end{align}
This implies that the family of laws of \( x^{n}(t) \) is weakly compact. Moreover, \eqref{t9} allows us to establish that the weak limit point of \( x^{n}(t) \) as \( n \to \infty \) constitutes a weak solution to the system \eqref{1r1}. This proof follows a standard approach; for detailed exposition, refer to \textbf{Appendix \uppercase\expandafter{\romannumeral1}} of \cite{ref3}. The  uniqueness in law of the weak solution of the finite dimensional system \eqref{1r1} is similar to the proof of Theorem 2.1 in \cite{ref68} and the  details can be found in the proof of uniqueness. Finally, by a standard argument (Lemma 3.1 in \cite{ref1}), \eqref{r3} has a  maximal local weak solution, and the weak solution is unique in law, i.e., any two such solutions with the same initial value have the same law.

Furthermore, let $\tau $ be the explosion time or life time. Then,  for any $t\in [0,\tau)$, by It$\hat{\text{o}} $'s formula and  \eqref{r6}, we have
\begin{align*}
		\left | x(t) \right | ^{2}\le  C_{T}\left \| \varphi \right \| _{h}^{2}+C\int_{0}^{t} (\underset{r\in[0,s]}{\sup}\left | x(r) \right | ^{2}+1)\text{d}s+2\int_{0}^{t} ( x(r))^{T} G (s,x(r))\text{d}W(s),
	\end{align*}
which implies that $\tau =\infty $ almost surely by the stochastic Gronwall lemma of \cite{ref000111}. Consequently, the unique maximal local weak solution $x(t)$ is non-explosive almost surely within any finite time. $\Box$
~\\

Building upon Theorem 3.1, we will conduct a detailed analysis of the existence, uniqueness of the solution, and other properties of system \eqref{r1} under condition \textbf{(H1)}$-$\textbf{(H5)}.  Next, we explore the properties of the solution to \eqref{r1} using the Galerkin-type approximation technique.
~\\
\\\textbf{Theorem 3.2.} \emph{Consider \eqref{r1}. Suppose that the assumptions \textbf{(H1)}$-$\textbf{(H5)} hold, then the following statements hold:
\begin{enumerate}[1)]
        \item  For any initial value $\varphi \in C^{h}_{U_{1}} $, system \eqref{r1} has  a probabilistic weak  solution  $u(t)$ and a segment process $u_{t}$ to \eqref{r3} with $u_{0}=\varphi$;
        \item The weak solution $u(t)$  is unique in law, i.e., any two such solutions with the same initial value have the same law;
		\item  Let $\{\varphi_{k}\}^{\infty}_{k=1}$ and $\varphi$ be a sequence with $\|\varphi_{k}-\varphi\|_{h}\to 0$. Let $u(t;\varphi)$ be the probabilistic weak  solution of \eqref{r1} with the initial value $\varphi$. Then we have
\[\lim_{k\to\infty} |\mathbb{E}\mathcal{L}(u(t;\varphi_{k})|_{[0,T]})-\mathbb{E}\mathcal{L}( u(t;\varphi)|_{[0,T]})| = 0,\]
where $\mathcal{L}$ is  an arbitrary real-valued, bounded continuous function defined on $C([0, T], U_{1})$.
	\end{enumerate}}

 We begin by applying the Galerkin projection technique to transform the system \eqref{r1} into a finite-dimensional system. Let's assume we have orthonormal bases  $\left \{ \eta  _{1},\eta  _{2},\eta  _{3} ,...\right \} \subset B$ for $U_{1}$ and $\left \{ \varrho_{1}, \varrho_{2}, \varrho_{3}, ..., \right \}$ for $U_{2}$. Selecting the first $k$ orthonormal bases from each set, we define the following operators:
\begin{align*}
\begin{split}
		\Theta _{1}^{k}:B^{*}\to U_{1}^{k}:=\text{span}\{\eta  _{1},\eta  _{2},...,\eta  _{k} \},\quad \Theta _{2}^{k}:U_{2}\to U_{2}^{k}:=\text{span}\{\varrho   _{1},\varrho   _{2},...,\varrho   _{k} \}.
\end{split}
	\end{align*}
For any $u\in B^{*}$ and $k\ge1$, we obtain $u^{k}=\Theta _{1}^{k}(u)=\sum_{i=1}^{k}~_{B ^{\ast}} \langle u, \eta  _{i}\rangle _{B}\eta  _{i}$ and  $W^{k}(t)=\Theta _{2}^{k}[W(t)]=\sum_{i=1}^{k}\left \langle W(t),\varrho   _{i} \right \rangle _{U_{2}}\varrho   _{i}$. Based on the above  approximation techniques, we will first analyze the finite-dimensional equation corresponding to system  \eqref{r1} as follows:
\begin{align}\label{r2}
		\begin{cases}
 \text{d}u^{k}(t)=\Theta _{1}^{k}[(A(t,u(t))+f(t,u_{t}))]\text{d}t+\Theta _{1}^{k}[g(t,u_{t})]\text{d}W^{k}(t), \\
u^{k}_{0} =\varphi^{k}\in C^{h,k}_{U_{1}},
\end{cases}
	\end{align}
where
\begin{eqnarray*}
		C^{h,k}_{U_{1}} =\left \{ \varphi \in C((-\infty  , 0]; U^{k}_{1}):\underset{\theta \to -\infty }{\lim }e^{h\theta }\varphi (\theta )=\mathcal{U}  \in   U^{k}_{1} \right \}.
	\end{eqnarray*}
Following Theorem 3.1, under assumptions \textbf{(H1)}, \textbf{(H3)}, \textbf{(H4)}, and \textbf{(H5)}, system \eqref{r2} admits a weak unique continuous solution $u^{k}(t)$. Next, we proceed to prove Theorem 3.2. Similarly, the proof is conducted in the following apriori estimates.
\\~
\\\textbf{Lemma 3.3.} \textbf{(Apriori estimates of the  solutions $u^{k}(t)$)} \emph{ Suppose that \textbf{(H1)}-\textbf{(H5)} hold. Then there exists a constant $L_{\alpha _{1},\alpha _{2},M,T}$, which is
independent of $k$, such that for any $k\ge 1$ and $T\ge 0$,
\begin{align*}
		&\mathbb{E} \underset{z\in[0,T]}{\sup}\left \| u^{k}(z) \right \|^{2}_{U_{1}} +\mathbb{E} \underset{z\in[0,T]}{\sup}\left \| u^{k}_{z} \right \|^{2}_{h}+\mathbb{E}\int _{0}^{T}[\left \| (A(s,u^{k}(s)) \right \|_{B^{*}}^{\frac{p}{p-1}}\nonumber
 \\&+\left \| f(s,u^{k}_{s}) \right \|^{2}_{U_{1}}
+\left \| g(s,u^{k}_{s}) \right \| ^{2}_{\mathscr{L}(U_{2},U_{1})}+\left \| u^{k}(s) \right \|_{B}^{p} ]\text{d}s
\\&\le  L_{\alpha _{1},\alpha _{2},M,T}(1+\left \| \varphi \right \|^{2}_{h }) \nonumber.
	\end{align*}
 \textbf{proof:}} By It\^o's formula for \eqref{r2}, \textbf{(H2)} and \textbf{(H3)}, we have
\begin{align}\label{p1}
		\left \| u^{k}(t) \right \|^{2}_{U_{1}} &=\left \| \varphi ^{k} \right \|^{2}_{h } +\int_{0}^{t}[2_{B ^{\ast}} \langle \Theta _{1}^{k}(A(s,u^{k}(s)), u^{k}(s)\rangle _{B}\nonumber
\\&~~~+2 \langle \Theta _{1}^{k}(f(s,u^{k}_{s}), u^{k}(s)\rangle _{U_{1}} +\left \|  \Theta _{1}^{k}[g(s,u^{k}_{s})]\Theta _{2}^{k}\right \|^{2}_{\mathscr{L}(U_{2},U_{1})} ]\text{d}s\nonumber
\\&~~~+2\int _{0}^{t}\langle u^{k}(s),\Theta _{1}^{k}[g(s,u^{k}_{s})]\text{d}W^{k}(s)\rangle_{U_{1}}
\\&\le\left \| \varphi ^{k} \right \|^{2}_{h }+\int_{0}^{t}[-\alpha _{1}\left \| u^{k}(s) \right \|_{B}^{p}+2\alpha _{2}\left \| u^{k}(s) \right \| _{U_{1}}^{2}+(2\alpha _{1}^{2}+2\alpha _{1}+1)\left \| u_{s}^{k} \right \| ^{2}_{h}+M]\text{d}s\nonumber
\\&~~~+2\int _{0}^{t}\langle u^{k}(s),\Theta _{1}^{k}[g(s,u^{k}_{s})]\text{d}W^{k}(s)\rangle_{U_{1}}.\nonumber
	\end{align}
 Note
\begin{align}\label{y1}
\left \| u^{k}_{s} \right \| _{h}\le\left \| \varphi  \right \|  _{h}+\underset{z\in[0,s]}{\sup}\left \| u^{k}(z) \right \|_{U_{1}}.
\end{align}
For arbitrary $N>0$, let $\tau ^{k}_{N}=\underset{t\ge0}{\inf} \{\left \| u^{k}(t)\right \|_{U_{1}}>N \}$. According  to  the Burkholder-Davis-Gundy  inequality, Young's inequality and \eqref{y1}, for arbitrary fixed time $T$ we get
\begin{align}\label{p2}
		&\mathbb{E} \underset{z\in[0,T\wedge \tau ^{k}_{N}]}{\sup}\left \| u^{k}(z) \right \|^{2}_{U_{1}} \nonumber
\\&\le L_{\alpha _{1},M,T}(1+\left \| \varphi ^{k} \right \|^{2}_{h })+L_{\alpha _{1},\alpha _{2}}\mathbb{E}\int_{0}^{T\wedge \tau ^{k}_{N}}\underset{z\in[0,s]}{\sup}\left \| u^{k}(z) \right \| _{U_{1}}^{2}\text{d}s \nonumber
\\&~~+6\mathbb{E}\int _{0}^{T\wedge \tau ^{k}_{N}} \left \| u^{k}(t) \right \| ^{2}_{U_{1}}\left \|\Theta _{1}^{k}[g(s,u^{k}_{s})]\right \| ^{2}_{\mathscr{L}(U_{2},U_{1})} \text{d}s
\\&\le L_{\alpha _{1},M,T}(1+\left \| \varphi ^{k} \right \|^{2}_{h })+L_{\alpha _{1},\alpha _{2}}\mathbb{E}\int_{0}^{T}\underset{z\in[0,s\wedge \tau ^{k}_{N}]}{\sup}\left \| u^{k}(z) \right \| _{U_{1}}^{2}\text{d}s  \nonumber
\\&~~~+\frac{1}{2}\mathbb{E} \underset{z\in[0,T\wedge \tau ^{k}_{N}]}{\sup}\left \| u^{k}(t) \right \|^{2}_{U_{1}} \nonumber.
	\end{align}
Then by Gronwall's lemma, we obtain
\begin{align}\label{p3}
		\mathbb{E} \underset{z\in[0,T\wedge \tau ^{k}_{N}]}{\sup}\left \| u^{k}(z) \right \|^{2}_{U_{1}} &\le L_{\alpha _{1},\alpha _{2},M,T}(1+\left \| \varphi \right \|^{2}_{h }),
	\end{align}
further,
\begin{align*}
\mathbb{E} \underset{z\in[0,T\wedge \tau ^{k}_{N}]}{\sup}\left \| u^{k}_{z} \right \|^{2}_{h} &\le L_{\alpha _{1},\alpha _{2},M,T}(1+\left \| \varphi \right \|^{2}_{h }).
\end{align*}
In addition, by taking expectations on both sides of \eqref{p1}, we derive
\begin{align}\label{p5}
		\mathbb{E}\int _{0}^{T\wedge \tau ^{k}_{N}}\left \| u^{k}(s) \right \|_{B}^{p}\text{d}s\le L_{\alpha _{1},\alpha _{2},M,T}(1+\left \| \varphi \right \|^{2}_{h }),
	\end{align}
and combining \textbf{(H2)}, \textbf{(H3)}, \eqref{p3} and \eqref{p5}, we obtain
\begin{align}\label{p6}
		&\mathbb{E}\int _{0}^{T\wedge \tau ^{k}_{N}}[\left \| (A(s,u^{k}(s)) \right \|_{B^{*}}^{\frac{p}{p-1}} +\left \| f(s,u^{k}_{s}) \right \|^{2}_{U_{1}}+\left \| g(s,u^{k}_{s} ) \right \| ^{2}_{\mathscr{L}(U_{2},U_{1})} ]\text{d}s\nonumber
\\&\le  L_{\alpha _{1},\alpha _{2},M,T}(1+\left \| \varphi \right \|^{2}_{h }).
	\end{align}
Taking $N\to\infty $, by the monotone convergence
theorem we have
\begin{align}\label{p7}
		&\mathbb{E} \underset{z\in[0,T]}{\sup}\left \| u^{k}(z) \right \|^{2}_{U_{1}} +\mathbb{E} \underset{z\in[0,T]}{\sup}\left \| u^{k}_{z} \right \|^{2}_{h}+\mathbb{E}\int _{0}^{T}[\left \| (A(s,u^{k}(s)) \right \|_{B^{*}}^{\frac{p}{p-1}}\nonumber
 \\&+\left \| f(s,u^{k}_{s}) \right \|^{2}_{U_{1}}
+\left \| g(s,u^{k}_{s}) \right \| ^{2}_{\mathscr{L}(U_{2},U_{1})}+\left \| u^{k}(s) \right \|_{B}^{p} ]\text{d}s
\\&\le  L_{\alpha _{1},\alpha _{2},M,T}(1+\left \| \varphi \right \|^{2}_{h }) \nonumber.
	\end{align}
$\Box$
~\\

Now we give the proof of Theorem 3.2.
~\\
\\\emph{\textbf{proof of Theorem 3.2:}}
\begin{enumerate}[1)]
        \item Due to the the reflexivity of $\left \| \cdot\right \|_{B}^{p}$ and \textbf{Lemma 3.3}, we may assume that there exist common subsequences $k_{n}$ such that for $n\to \infty  $:
\begin{enumerate}[\textbf{(1)}]
		\item $u^{k_{n}}(t)\to u(t)$ in $ L ^{2}([0,T]\times \Omega ,U_{1})$ and  weakly in $ L ^{p}([0,T]\times \Omega ,B)$;
\end{enumerate}
\begin{enumerate}[\textbf{(2)}]
        \item $A(t,u^{k_{n}}(t) )\to A^{*}(t)$ weakly in $[ L ^{p}([0,T]\times \Omega ,B)]^{*}$;
	\end{enumerate}
\begin{enumerate}[\textbf{(3)}]
        \item $f(t,u^{k_{n}}_{t})\to f^{*}(t)$ weakly in $ L ^{2}([0,T]\times \Omega ,U_{1})$;
	\end{enumerate}
\begin{enumerate}[\textbf{(4)}]
        \item $g(t,u^{k_{n}}_{t})\to g^{*}(t)$ weakly in $ L ^{2}([0,T]\times \Omega ,\mathscr{L}(U_{2},U_{1}))$.
	\end{enumerate}
Then for any $v \in B  $ and $t\in[0,T]$,
\begin{align*}
		&\mathbb{E} \int_{0}^{t}  [_{B ^{\ast}} \langle u(s), v \rangle _{B}]\text{d}s
\\&=\lim_{n \to \infty} \mathbb{E} \int_{0}^{t}  [_{B ^{\ast}} \langle u^{k_{n}}(s), v \rangle _{B}]\text{d}s
\\&=\lim_{n \to \infty}\mathbb{E}\int_{0}^{t}[_{B ^{\ast}} \langle \varphi^{k_{n}} (0), v \rangle _{B}+\int_{0}^{s}  (_{B ^{\ast}} \langle A(z,u^{k_{n}}(z)), v \rangle _{B})\text{d}z
\\&~~~+\int_{0}^{s}\left \langle f(z,u^{k_{n}}_{z}),v \right \rangle _{U_{1}}\text{d}z
+\int _{0}^{s}\langle v ,g(z,u^{k_{n}}_{z})\text{d}W(z)\rangle_{U_{1}}]\text{d}s
\\&=\mathbb{E}\int_{0}^{t}[_{B ^{\ast}} \langle \varphi (0), v \rangle _{B}+\int_{0}^{s}  (_{B ^{\ast}} \langle A^{*}(z), v \rangle _{B})\text{d}z
\\&~~~+\int_{0}^{s}\left \langle f^{*}(z),v \right \rangle _{U_{1}}\text{d}z
+\int _{0}^{s}\langle v,g^{*}(z)\text{d}W(z)\rangle_{U_{1}}]\text{d}s.
	\end{align*}
Therefore, for any $t\in [0,T]$,
\begin{align}\label{y3}
		\begin{cases}
 u(t)=\varphi (0)+\int_{0}^{t}A^{*}(s)\text{d}s+\int_{0}^{t}f^{*}(s)\text{d}s+\int_{0}^{t}g^{*}(s)\text{d}W(s), \quad \text{d}t\times \mathbb{P}-\text{a.e.}, \\
u_{0} =\varphi\in C^{h}_{U_{1}}.
\end{cases}
	\end{align}
In addition, let  $$\psi (t) \in  L ^{p}((- \infty ,T]\times \Omega ,B)\cap  L ^{2}((- \infty ,T]\times \Omega ,U_{1}),$$
 where  $\psi (\theta)=\varphi (\theta ), \theta\in(- \infty ,0]$. 
%
Then the It\^o's formula yields
\begin{align}\label{p9}
		 &\mathbb{E}\left \| u^{k_{n}}(t) \right \|_{U_{1}}^{2}  -\left \|\varphi (0) \right \|_{U_{1}}^{2}   \nonumber
\\&=\mathbb{E} \int_{0}^{t}[2_{B ^{\ast}} \langle A(s,u^{k_{n}}(s)), u^{k_{n}}(s) \rangle _{B}+2\langle f(s,u^{k_{n}}_{s}), u^{k_{n}}(s) \rangle _{U_{1}}  +\left \| g(s,u^{k_{n}}_{s}) \right \|_{\mathscr{L}(U_{2},U_{1})}^{2}]\text{d}s  \nonumber
\\&\le \mathbb{E} \int_{0}^{t}[2_{B ^{\ast}} \langle A(s,u^{k_{n}}(s))-A(s,\psi (s)), u^{k_{n}}(s)-\psi (s) \rangle _{B}
\\&~~~+2\langle f(s,u^{k_{n}}_{s})-f(s,\psi _{s}), u^{k_{n}}(s)-\psi (s) \rangle _{U_{1}}
+\left \| g(s,u^{k_{n}}_{s})-g(s,\psi _{s}) \right \|_{\mathscr{L}(U_{2},U_{1})}^{2} \nonumber
\\&~~~+2_{B ^{\ast}} \langle A(s,\psi (s)), u^{k_{n}}(s) \rangle _{B}+2\langle f(s,\psi _{s}), u^{k_{n}}(s) \rangle _{U_{1}}   \nonumber
\\&~~~+2_{B ^{\ast}} \langle A(s,u^{k_{n}}(s) )-A(s,\psi (s)), \psi (s) \rangle _{B}
+2\langle f(s,u^{k_{n}}_{s})-f(s,\psi _{s}), \psi (s) \rangle _{U_{1}}  \nonumber
\\&~~~+2\langle g(s,u^{k_{n}}_{s}),g(s,\psi _{s}) \rangle _{\mathscr{L}(U_{2},U_{1})}-\left \|g(s,\psi _{s})  \right \|_{\mathscr{L}(U_{2},U_{1})}^{2}]\text{d}s,   \nonumber
	\end{align}
where $\psi_{s}=:\{\psi(s+\theta  ), \theta  \in(- \infty  , 0]\}$.  Without loss of generality,  there exists a sufficiently large constant $M_{T}$ such that $\mathbb{E} \left \| \psi (z) \right \|^{2}_{U_{1}}\le M_{T}$.

~~~By \textbf{(H4)}, H\"{o}lder inequality and Jensen's inequality,  we have
\begin{align*}
		  \Xi  (t):&=\mathbb{E} \int_{0}^{t}[2_{B ^{\ast}} \langle A(s,u^{k_{n}}(s))-A(s,\psi (s)), u^{k_{n}}(s)-\psi (s) \rangle _{B}   \nonumber
\\&~~~+2\langle f(s,u^{k_{n}}_{s})-f(s,\psi _{s}), u^{k_{n}}(s)-\psi (s) \rangle _{U_{1}}
+\left \| g(s,u^{k_{n}}_{s})-g(s,\psi _{s}) \right \|_{\mathscr{L}(U_{2},U_{1})}^{2}]\text{d}s
\\&\le\mathbb{E} \int_{0}^{t}[\alpha_{1}\left \| u^{k_{n}}_{s}-\psi _{s} \right \| _{h}^{2  }+L_{M}(\left \| u^{k_{n}}_{s}-\psi _{s} \right \| _{h}^{\beta +1}+\left \| u^{k_{n}}_{s}-\psi _{s} \right \| _{h}^{2\gamma  })]\text{d}s
\\&\le \alpha_{1}\mathbb{E} \int_{0}^{t}\left \| u^{k_{n}}_{s}-\psi _{s} \right \| _{h}^{2  }\text{d}s+L_{M,T}[(\mathbb{E} \int_{0}^{t}\left \| u^{k_{n}}_{s}-\psi _{s} \right \|_{h}^{2}\text{d}s)^{\frac{\beta +1}{2}}+(\mathbb{E} \int_{0}^{t}\left \| u^{k_{n}}_{s}-\psi _{s} \right \| _{h}^{2  }\text{d}s)^{\gamma}].
	\end{align*}
Hence, for  any given nonnegative function $\kappa   \in L^{\infty }([0,T],\mathbb{R} )$ and letting $n \to \infty $, it follows from \eqref{p9} that
\begin{align}\label{p10}
		 &\mathbb{E}\int_{0}^{T}\kappa  (t) [\left \| u(t) \right \|_{U_{1}}^{2}  -\left \|\varphi (0) \right \|_{U_{1}}^{2} ]\text{d}t  \nonumber
\\&\le\mathbb{E} \int_{0}^{T}\kappa  (t)[\int_{0}^{t}[2_{B ^{\ast}} \langle A(s,\psi (s)), u(s) \rangle _{B}+2\langle f(s,\psi _{s}), u(s) \rangle _{U_{1}} \nonumber
\\&~~~+2_{B ^{\ast}} \langle A^{*}(s)-A(s,\psi (s)), \psi (s) \rangle _{B}
+2\langle f^{*}(s )-f(s,\psi _{s}), \psi (s) \rangle _{U_{1}}
\\&~~~+2\langle g^{*}(s ), g(s,\psi _{s}) \rangle _{\mathscr{L}(U_{2},U_{1})}-\left \|g(s,\psi _{s})  \right \|_{\mathscr{L}(U_{2},U_{1})}^{2} ]\text{d}s\text{d}t \nonumber
 +\int_{0}^{T}\kappa  (t)[\alpha_{1}\mathbb{E} \int_{0}^{t}\left \| u_{s}-\psi _{s} \right \| _{h}^{2  }\text{d}s
 \\&~~~+L_{M,T}[(\mathbb{E} \int_{0}^{t}\left \| u_{s}-\psi _{s} \right \|_{h}^{2}\text{d}s)^{\frac{\beta +1}{2}}+(\mathbb{E} \int_{0}^{t}\left \| u_{s}-\psi _{s} \right \| _{h}^{2  }\text{d}s)^{\gamma}]]\text{d}t.   \nonumber
	\end{align}
By It\^o's formula to $\left \| u(t) \right \|_{U_{1}}^{2}  -\left \|\varphi (0) \right \|_{U_{1}}^{2}$, we have
\begin{align}\label{p12}
		 &\mathbb{E} \int_{0}^{T}\kappa  (t)\int_{0}^{t}[2_{B ^{\ast}} \langle A^{*}(s )-A(s,\psi (s)), u(s)-\psi (s) \rangle _{B}  \nonumber
\\&+2\langle  f^{*}(s )-f(s,\psi _{s}), u(s)-\psi (s) \rangle _{U_{1}}+\left \|g^{*}(s )-g(s,\psi _{s})  \right \|_{\mathscr{L}(U_{2},U_{1})}^{2}
]\text{d}s\text{d}t
 \\&-\int_{0}^{T}\kappa  (t)[\alpha_{1}\mathbb{E} \int_{0}^{t}\left \| u_{s}-\psi _{s} \right \| _{h}^{2  }\text{d}s
+L_{M,T}[(\mathbb{E} \int_{0}^{t}\left \| u_{s}-\psi _{s} \right \|_{h}^{2}\text{d}s)^{\frac{\beta +1}{2}}+(\mathbb{E} \int_{0}^{t}\left \| u_{s}-\psi _{s} \right \| _{h}^{2  }\text{d}s)^{\gamma}]]\text{d}t\nonumber
\\&~~~ \le 0.   \nonumber
	\end{align}
 First, taking $\psi (t)=u(t)$ for $t\ge 0$ implies that $g^{*}(t )=g(t,u _{t})$. In addition, let $\psi =u-\varepsilon \kappa  ^{*}y$ where $\kappa  ^{*} \in L^{\infty }((-\infty ,T],\mathbb{R} )$, $\kappa  ^{*}(\theta)=0$ for $\theta\in (-\infty ,0]$, $\varepsilon >0$ and $y\in B$, then we have
\begin{align*}
		 &\mathbb{E} \int_{0}^{T}\kappa  (t)\int_{0}^{t}[2_{B ^{\ast}} \langle A^{*}(s )-A(s,u(s)-\varepsilon \kappa  ^{*}(s)y ), \varepsilon \kappa  ^{*}(s)y \rangle _{B} \nonumber
\\&+2\langle  f^{*}(s )-f(s,u_{s}-\varepsilon \kappa  ^{*}_{s}y), \varepsilon \kappa  ^{*}(s)y \rangle _{U_{1}} ]\text{d}s\text{d}t-\int_{0}^{T}\kappa  (t)[\alpha_{1}\mathbb{E} \int_{0}^{t}\left \| \varepsilon \kappa  ^{*}_{s}y \right \| _{h}^{2  }\text{d}s
\\&+L_{M,T}[(\mathbb{E} \int_{0}^{t}\left \| \varepsilon \kappa  ^{*}_{s}y \right \|_{h}^{2}\text{d}s)^{\frac{\beta +1}{2}}+(\mathbb{E} \int_{0}^{t}\left \| \varepsilon \kappa  ^{*}_{s}y \right \| _{h}^{2  }\text{d}s)^{\gamma}]]\text{d}t\nonumber
\\&\le0.
	\end{align*}
Similarly, the converse follows by taking $ \kappa  ^{*}(s)=- \kappa  ^{*}(s)$, and finally  according to Lebesgue's dominated convergence theorem, we obtain that when $\varepsilon\to 0$,
\begin{align*}
		 &\mathbb{E} \int_{0}^{T}\kappa  (t)\int_{0}^{t}[2_{B ^{\ast}} \langle A^{*}(s )-A(s,u(s)), \kappa  ^{*}(s)y \rangle _{B} \nonumber
\\&+2\langle  f^{*}(s )-f(s,u_{s}), \kappa  ^{*}(s)y \rangle _{U_{1}} ]\text{d}s\text{d}t
\\&=0,
	\end{align*}
 which concludes $A^{*}(t)=A(t,u(t  ))$ and  $ f^{*}(t)=f(t,u_{t  })$, i.e., it suffices to prove that
\begin{align*}
A^{*}=A(\cdot,u(\cdot  )), \quad f^{*}=f(\cdot, u_{(\cdot)  }),\quad g^{*}=g(\cdot,u_{(\cdot)  }), \quad \text{d}t\times \mathbb{P}-\text{a.e.}
	\end{align*}
By \eqref{y3}, this completes the existence proof, i.e.,
\begin{align*}
\begin{cases}
 u(t)=\varphi (0)+\int_{0}^{t}A(s,u(s  ))\text{d}s+\int_{0}^{t}f(s,u_{s  })\text{d}s+\int_{0}^{t}g(s,u_{s  })\text{d}W(s), \quad \text{d}t\times \mathbb{P}-\text{a.e.}, \\
u_{0} =\varphi\in C^{h}_{U_{1}}.
\end{cases}
	\end{align*}
 \item To establish the weak uniqueness of  $u(t)$ for system \eqref{r1}, we  employ the  generalized coupling approach. The proof is inspired by \cite{ref68}.
 We choose a family of finite-dimensional projections $\{\Lambda _{n}\}_{n\ge1 }$ in $C^{h}_{U_{1}}$, which possesses the following property:
\begin{align*}
		\Lambda _{n}\to I(n\to\infty )\quad\text{and}\quad \Lambda  _{n}[\varphi (0)]=\varphi (0),
	\end{align*}
where $I$ are the identity transformation on $C^{h}_{U_{1}}$. Let $$f^{n,\delta}(t,\varphi)=f(t,\Lambda_{n}(\varphi) ), \quad  g^{n,\delta}(t,\varphi)=g(t,\Lambda_{n}(\varphi) ),$$ then we convolve the $(f^{n,\delta},g^{n,\delta})$ with the finite-dimensional approximation $\delta $-function to obtain the $(f^{n},g^{n})$:
\begin{enumerate}[(\textbf{I})]
		\item $f^{n}\to f$, $g^{n}\to g$ as $n\to \infty $  uniformly on each compact subset of $ C^{h}_{U_{1}}$;
\end{enumerate}
\begin{enumerate}[(\textbf{II})]
        \item $f^{n},g^{n}$ satisfy the conditions \textbf{(H1)}-\textbf{(H3)}, that is,  the constants are independent of $n$;
	\end{enumerate}
\begin{enumerate}[(\textbf{III})]
        \item the functions $f^{n},g^{n}$ are Lipschitz continuous on each bounded subset of $C^{h}_{U_{1}}$.
	\end{enumerate}
Now consider the following SFPDEs with infinite delay:
\begin{align}
    \begin{cases}
        \text{d}u(t) = \left(A(t,u(t))+f(t,u_{t})\right)\text{d}t+g(t,u_{t})\text{d}W(t), \\
         u_{0} =\varphi\in C^{h}_{U_{1}},
    \end{cases} \label{r1r}
\end{align}
\begin{align}
    \begin{cases}
        \text{d}u^{n}(t) = \left((A(t,u^{n}(t)) + f^{n}(t,u^{n}_{t})\right) \text{d}t + g^{n}(t,u^{n}_{t}) \text{d}W(t), \\
        u_{0} =\varphi\in C^{h}_{U_{1}},
    \end{cases} \label{r2r}
\end{align}
and
\begin{align}
    \begin{cases}
        \text{d}u^{n,\zeta}(t) =  \left(A(t,u^{n,\zeta}(t)) + f^{n}(t,u^{n,\zeta}_{t})+ \zeta(u(t) - u^{n,\zeta}(t))\mathbf{1} _{t\le\tau}\right) \text{d}t
         + g^{n}(t,u^{n,\zeta}_{t}) \text{d}W(t)
       \text{d}t, \\
        u_{0} =\varphi\in C^{h}_{U_{1}},
    \end{cases} \label{r3r}
\end{align}
where  the constant $\zeta$ and the stopping time $\tau$  will be
determined later. Given any compact subset  $\mathcal{C} \subset C^{h}_{U_{1}}$,  denote
\[
\tau_{1}=\inf \{t\ge 0;u_{t}\notin \mathcal{C}\}, \quad
 \Re_{f}=\sup_{t\ge0}\sup_{\varphi\in \mathcal{C}}\| f^{n}(t,\varphi)-f(t,\varphi)\|_{U_{1}},
\]
\[
\Re_{g}=\sup_{t\ge0}\sup_{\varphi\in \mathcal{C}}\| g^{n}(t,\varphi)-g(t,\varphi)\|_{\mathscr{L}(U_{2},U_{1})},\quad \Re=\max\{\Re_{f}^{\frac{1}{\beta}} ,\Re_{g}^{\frac{1}{\gamma}}  \},
\]
\[\tau_{2}=\inf\{t\ge 0; \|u(t)-u^{n,\zeta}(t)\|_{U_{1}}>2\Re\}.\]
Let \( \tau= \tau_{1} \wedge \tau_2 \)  and $\zeta=\Re^{\ell-1} > 0$ be a constant with  $\ell\in (0, \beta\wedge(2\gamma-1))$.
%
%

~~~Now, applying Girsanov's theorem, we establish relationships between solutions $u^{n}(t)$ and $u^{n,\zeta}(t)$. Let $d_{\text{TV}}$ denote total variation distance and $KL(\cdot\|\cdot)$ the Kullback-Leibler divergence on $C([0,T];U_1)$. The strong monotonicity of coefficients $f^n$ and $g^n$ ensures the existence of strong solutions $u^n(t)$ and $u^{n,\zeta}(t)$. The controlled system \eqref{r2r} relates to \eqref{r3r} through a shifted Wiener process:
\[
\text{d}W^{*}(t) = \text{d}W(t) + \mathcal{W}(t) \text{d}t, \quad \mathcal{W} (t) = g^{n}(t,u^{n,\zeta}_{t})^{-1} \zeta(u(t) - u^{n,\zeta}(t))\mathbf{1} _{t\le\tau\wedge T} \text{d}t.
\]
This gives
\[
d_{\text{TV}} \left( \text{Law}(u^{n}(t)| _{[0,T]}), \text{Law}(u^{n,\zeta}(t) | _{[0,T]}) \right) \leq d_{\text{TV}} \left( \text{Law}(W(t) | _{[0,T]}), \text{Law}(W^{*} (t)| _{[0,T]}) \right).
\]

Since $g^{-1}$ is  bounded, by recalling the definition of the stopping time $\tau$ and   control term $\zeta$, applying  Pinsker's inequality, we deduce that
\begin{align*}
&2[d_{TV}(\text{Law}(W^{n}(t) | _{[0,T]}), \text{Law}(W^{*}(t)| _{[0,T]}))]^{2}\\&\le KL (\text{Law}(W^{n}(t)| _{[0,T]}) \| \text{Law}(W^{*}(t)|_{ [0,T]}))
\\&\le \mathbb{E} \left[ \int_0^\infty \left\| \mathcal{W}(t) \right\|_{U_{1}}^2dt \right]
\\&\leq C_{T} \Re^{2\ell}.
\end{align*}
The total variation distance satisfies the following bound:
\begin{align}\label{r10r}
d_{\text{TV}} \left( \text{Law}(u^{n}(t)| _{[0,T]}), \text{Law}(u^{n,\zeta}(t) | _{[0,T]}) \right)  \leq C _{T} \Re^{\ell}.
\end{align}

 ~~~In the following, we estimate $u(t) - u^{n,\zeta}(t)$ via pathwise arguments.  For any \( t \leq \tau \), we have $u_{t}\in \mathcal{C}$. Applying Itô's formula, \textbf{(H4)} and \eqref{z1},  we obtain
\begin{align}\label{r5r}
&\left\| u(t) - u^{n,\zeta}(t) \right\|_{U_{1}}^2\nonumber
\\&=  \int_0^t [2_{B ^{\ast}} \langle A(s,u(s))-A(s,u^{n,\zeta}(s)), u(s) - u^{n,\zeta}(s)\rangle _{B}+ 2\langle f(s,u_{s}) - f^{n}(s,u^{n,\zeta}_{s}), u(s) - u^{n,\zeta}(s) \rangle_{U_{1}}\nonumber
\\&~~~+ \left\| g(s,u_{s}) - g^{n}(s,u^{n,\zeta}_{s}) \right\|_{\mathscr{L}(U_{2},U_{1})}^2- 2 \zeta\|u(s) - u^{n,\zeta}(s)\|_{U_{1}}^{2}\mathbf{1} _{t\le\tau}]\text{d}s\nonumber
\\&~~~ + 2\int _{0}^{t}\langle u(s) - u^{n,\zeta}(s),[g(s,u_{s}) - g^{n}(s,u^{n,\zeta}_{s})]\text{d}W(s)\rangle_{U_{1}}\nonumber
\\&\le  \int_0^t [2_{B ^{\ast}} \langle A(s,u(s))-A(s,u^{n,\zeta}(s)), u(s) - u^{n,\zeta}(s)\rangle _{B}+ 2\langle f(s,u_{s}) - f^{n}(s,u_{s}), u(s) - u^{n,\zeta}(s) \rangle_{U_{1}}\nonumber
\\&~~~+ \left\| g(s,u_{s}) - g^{n}(s,u_{s})  \right\|_{\mathscr{L}(U_{2},U_{1})}^2+ 2\langle f^{n}(s,u_{s}) - f^{n}(s,u^{n,\zeta}_{s}), u(s) - u^{n,\zeta}(s) \rangle_{U_{1}}
\\&~~~+\left\| g^{n}(s,u_{s}) - g^{n}(s,u^{n,\zeta}_{s})  \right\|_{\mathscr{L}(U_{2},U_{1})}^2- 2 \zeta\|u(s) - u^{n,\zeta}(s)\|^{2}_{U_{1}}\mathbf{1} _{t\le\tau}]\text{d}s\nonumber
\\&~~~ + 2\int _{0}^{t}\langle u(s) - u^{n,\zeta}(s),[g(s,u_{s}) - g^{n}(s,u_{s}) ]\text{d}W(s)\rangle_{U_{1}}\nonumber
\\&~~~ + 2\int _{0}^{t}\langle u(s) - u^{n,\zeta}(s),[g^{n}(s,u_{s}) - g^{n}(s,u^{n,\zeta}_{s}) ]\text{d}W(s)\rangle_{U_{1}}\nonumber
\\&\le C\int _{0}^{t}[(\Re^{2}+\Re_{f}\Re+\Re_{g}^{2}+\Re^{\beta}\|u(s) - u^{n,\zeta}(s)\|_{U_{1}}+\Re^{2\gamma})-2 \Re^{\ell-1}\|u(s) - u^{n,\zeta}(s)\|^{2}_{U_{1}}]\text{d}t\nonumber
\\&~~~+C\int _{0}^{t}(\Re_{g}\Re+\Re^{\gamma+1})\text{d}W(s).\nonumber
\end{align}
By Young's inequality,
\begin{align*}
\Re_{f}\|u(s) - u^{n,\zeta}(s)\|_{U_{1}}\le \Re^{\beta}\|u(s) - u^{n,\zeta}(s)\|_{U_{1}}
\le \frac{1}{2}[\Re^{\beta+1}+\Re^{\beta-1}\|u(s) - u^{n,\zeta}(s)\|^{2}_{U_{1}}].
\end{align*}
On the other hand, since \( \ell < \beta \), there exists \( \Re^{*} > 0 \) such that for any \( t \in[0, \tau] \) and $\Re\in (0,\Re^{*}]$,
\begin{align}\label{r9r}
\left\| u(t) - u^{n,\zeta}(t) \right\|_{U_{1}}^2\nonumber
\le C\int _{0}^{t}[-\Re^{\ell-1}\|u(t) - u^{n,\zeta}(t)\|^{2}_{U_{1}}
+\Re^{\beta+1}+\Re^{2\gamma}+\Re^{2}]\text{d}t+ C\int _{0}^{t}\Re^{\gamma+1}\text{d}W(s).
\end{align}
Now we apply Lemma 2.1 with a fixed \( T > 0 \), then we obtain
\[
\mathbb{P}\left( \sup_{t \leq \tau \wedge T}( \left\|u(t) - u^{n,\zeta}(t) \right\|_{U_{1}}^2 )  \geq C \left(\Re^{\beta+1}+\Re^{2\gamma}+\Re^{2} \right)\Re^{(1-\ell)}+C \Re^{\gamma+1}\Re^{(1-\ell)\delta}R \right)\leq C_{1}e^{-C_{2}R^{2}},
\]
for   $R\geq 0$ and every $\delta\in (0,\frac{1}{2})$.
Recall that $\ell\in (0, \beta\wedge(2\gamma-1))$; hence, there exists \( \xi> 0 \) such that for for sufficiently small \( \Re^{*} > 0 \) and any $\Re\in (0,\Re^{*}]$, we obtain that
\begin{align*}
C[\Re^{\beta+2-\ell}+\Re^{2\gamma+1-\ell}+\Re^{3-\ell}]
\le\Re^{2+2\xi}.
\end{align*}
Furthermore, it follows that there exists $\delta_{0}\in (0,\frac{1}{2})$ such that $\gamma+\delta_{0}(1-\ell)>1$. Moreover, let $\xi\in (0,\frac{\gamma+\delta_{0}(1-\ell)-1}{3})$. Then we have
\[
B^{\frac{1}{2}}\lambda^{-\delta_{0}} \leq \varepsilon \Re^{2 + 3\xi}.
\]
Let  \( R = \Re^{-\ell} \),  and denote
\[
\mathbb{A} = \left\{ \omega : \sup_{t \leq \tau \wedge T}   \left\|u(t) - u^{n,\zeta}(t) \right\|_{U_{1}}  \geq\Re^{1 + \xi} \right\}.
\]
 Without loss of generality, taking \( \Re^{*} \leq 1 \), we have for \( \Re \in (0,  \Re^{*}] \), $\Re^{\xi}\leq 1$. Then, on the set \( \Omega\setminus\mathbb{A}\), we have
\[
 \left\|u(t) - u^{n,\zeta}(t) \right\|^{2}_{U_{1}}  \leq\Re^{2 + 2\xi}\leq4\Re^{2}, \quad t \leq \tau \wedge T.
\]
Recall the definition of the stopping time \( \tau \). Given that  \( u(t),u^{n,\zeta}(t) \) exhibit continuous trajectories in the function space, it follows that on \( \Omega\setminus\mathbb{A} \), the stopping time satisfies \( \tau_{2} > \tau \wedge T \). Consequently, we establish that
\(
\tau \wedge T = \tau_{1} \wedge T.
\)
We thus derive the following result:
\[
\mathbb{P }\left(\sup_{t \in [0, \tau_{1}\wedge T]}   \left\|u(t) - u^{n,\zeta}(t) \right\|_{U_{1}}  \geq\Re^{1 + \xi} \right) \leq C_1 e^{-C_2 \Re^{-2\xi} }.
\]
This implies that for any $\varepsilon > 0$, as $n$ gets large and $\Re \to 0$, we obtain
\begin{align}\label{r11r}
 \mathbb{P}\left( \sup_{t \in [0, \tau_{1}\wedge T]} \left\|u(t) - u^{n,\zeta}(t) \right\|_{U_{1}}  > \varepsilon \right) \to 0, \quad n \to \infty.
\end{align}

~~~We can now complete  the proof of uniqueness in law. Consider $\mathcal{L}$ to be a real-valued, bounded continuous function defined on $C([0, T], U_{1})$. By \eqref{r10r}, we have
\begin{align}\label{b1b}
\mathbb{E}\mathcal{L}[u^{n}(t)] - \mathbb{E}\mathcal{L}[u^{n,\zeta}(t)] \to 0, \quad n \to \infty,
\end{align}
for any $t\in [0,T]$. In addition, when $ T\leq\tau_{1} $, \eqref{r11r} implies that for any  $\omega\in \Omega \setminus \mathbb{A}$, we have
\begin{align}\label{b2b}
u^{n,\zeta}(t) \to u(t) ,\quad  n \to \infty,
\end{align}
 in probability in $C([0, T], U_{1})$. Note that
\begin{align}\label{b3b}
 |\mathbb{E}\mathcal{L}(u^{n,\zeta}(t)) - \mathbb{E}\mathcal{L}(u(t))| \leq 2 [\sup_x |\mathcal{L}(x)| [\mathbb{P}(\tau_{1}\leq T)+\mathbb{P}(\tau_{1}\geq T)],
\end{align}
for any $t\in [0,T]$. Since $\lim_{n \to \infty} \mathbb{P}(\mathbb{A}) = 0$, by \eqref{b1b}-\eqref{b3b}, we obtain
\begin{align}\label{r13r}
\limsup_{n \to \infty} |\mathbb{E}\mathcal{L}(u^{n}(t)| _{[0,T]}) - \mathbb{E}\mathcal{L}(u(t)|_{[0,T]})| \leq 2 \sup_x |\mathcal{L}(x)| \mathbb{P}(\tau_{1} \leq T).
\end{align}
Similar to \textbf{Lemma 3.3}, we obtain that for any $p > 1$, there exists a constant $C>0$ such that
\begin{equation}\label{eq:condition}
\mathbb{E}\left[ \sup_{t \in [0,T]} \|u(t)\|_{U_{1}}^2 \right] + \mathbb{E}\left[ \int_0^T \|u(s)\|_{B}^{p} \, ds \right] \leq C .
\end{equation}
For any \( \varepsilon > 0 \), let \( M = \left(\frac{2C}{\varepsilon}\right)^{1/p} \) and  by Chebyshev's inequality, we have
\begin{equation}\label{eq:integral}
\mathbb{P}\left( \int_0^T \|u(s)\|_{B}^p ds \leq M^p \right) \geq 1 - \frac{\varepsilon}{2}.
\end{equation}
Given that \( \mathbb{E}[\sup_{t \in [0,T]} \|u(t)\|_{U_{1}}^2] \leq C \), assuming condition \textbf{(H3)} holds,  it follows from a standard procedure using the Kolmogorov continuity theorem that there exists a parameter  \( \iota\in (0,1/2) \) and a random variable \( \mathcal{V} \) satisfying \( \mathbb{E}\mathcal{V}^2 \leq  C \) such that for any $t,s\in[0,T]$,
\[
\|u(t) \|_{U_{1}}\le \mathcal{V},\quad \|u(t) - u(s)\|_{U_{1}} \leq \mathcal{V}|t-s|^\iota \quad \text{a.s.}
\]
Applying Chebyshev's inequality to \( \mathcal{V}  \) yields:
\[
\mathbb{P}\left( \mathcal{V}  > K \right) \leq \frac{\mathbb{E}\mathcal{V} ^2}{K^2} \leq \frac{ C}{K^2}, \quad \text{for } K > 0.
\]
Choosing \( K = \sqrt{\frac{2C}{\varepsilon}} \), we have
\begin{equation}\label{eq:holder}
\mathbb{P}\left( \mathcal{V}  \leq K \right) \geq 1 - \frac{\varepsilon}{2}.
\end{equation}
Then define the space $\mathbb{S} = L^p([0,T], B) \cap C^{0,\iota}([0,T], U_{1})$, and
using  the Aubin-Lions lemma and compact embedding \( B \subset U_{1} \), we conclude that bounded sets in \( \mathbb{S}  \) are compact in \( C([0,T], U_{1}) \). Define
\[
 \mathbb{E}= \left\{ u \in C([0,T], U_{1}) \,\Big|\, \|u\|_{C^{0,\iota}} \leq K, \, \int_0^T \|u(s)\|_{B}^p ds \leq M^p \right\},
\]
then \( \mathbb{E} \) is compact in \( C([0,T], U_{1}) \). For any \( t \in [0,T] \), consider the projection map \( \varpi_t: C([0,T], U_{1}) \to U_{1} \), \( \varpi_t(u) = u(t) \). Since \( \mathbb{E} \) is compact and \(\varpi_t \) is continuous, the image $
\mathbb{U}= \varpi_t(\mathcal{E}) = \{ u(t) \in U_{1} \mid u \in \mathbb{E}\}
$
is compact in \( U_{1} \). Combining \eqref{eq:integral} and \eqref{eq:holder}, we obtain
\[
\mathbb{P}\left( u \in \mathbb{E} \right) \geq 1 - \frac{\varepsilon}{2} - \frac{\varepsilon}{2} = 1 - \varepsilon.
\]
Therefore
\begin{align}\label{r15r}
\mathbb{P}\left( u(t) \in \mathbb{U}\right) \geq 1 - \varepsilon, \quad \forall t \in [0,T].
\end{align}
In the following discussion, we aim to demonstrate that the corresponding weak convergence is uniform with respect to $\varphi\in \mathcal{C}_{0}$ chosen from a compact set $\mathcal{C}_{0}\subseteq C^{h}_{U_{1}}$. To facilitate this analysis, we  define $\mathcal{C}$
  as the set of $\phi\in \mathcal{C}$ such that there exists an $r\in [0,T]$,
  \[\left\{\begin{matrix}
 \phi (\theta )=\varphi (\theta ), & \theta \in (-\infty ,-r],\\
 \phi (\theta )\in\mathbb{U},  & \theta \in [-r ,0].
\end{matrix}\right.\]
Then $\mathcal{C}$ is a compact subset of $C^{h}_{U_{1}}$, and it follows from \eqref{r15r} that for any $T > 0$,
\begin{align}\label{i15i}
\mathbb{P}\left( u_{t} \in \mathcal{C}\right) \geq 1 - \varepsilon, \quad \forall t \in [0,T].
\end{align}
This implies that for any $\varepsilon>0$, there exists a compact subset $\mathcal{C}\subseteq C^{h}_{U_{1}}$ such that with probability at least $1-\varepsilon$, the process $u_{t}$ remains in $\mathcal{C}$ for all $t\in[0, T]$. Consequently, we may select a compact set $\mathcal{C}$ ensuring that $\mathbb{P}(\tau_{1} \leq T)$ can be made arbitrarily small. Letting $\varepsilon\to 0$ and employing estimates \eqref{r13r} and \eqref{i15i}, we conclude that
\[\mathbb{E}\mathcal{L}(u^{n}(t)| _{[0,T]}) - \mathbb{E}\mathcal{L}(u(t)|_{[0,T]})\to 0,\quad n\to \infty,\]
thus proving weak uniqueness. This conclusion follows from the key observation that any weak solution to \eqref{r1r} can be uniquely identified  on arbitrary finite intervals $[0,T]$ as the weak limit of a sequence of solutions to  \eqref{r2r}. Consequently, the probability distribution of $u$ is uniquely determined.
\item To prove the continuity of the solution with respect to the initial data, we begin by considering the following system:
\begin{align*}
        \text{d}u^{n}(t) = \left((A(t,u^{n}(t)) + f^{n}(t,u^{n}_{t})\right) \text{d}t + g^{n}(t,u^{n}_{t}) \text{d}W(t).
\end{align*}
Similar to \textbf{Lemma 3.3}, we obtain that for any $p >1$, 
\begin{align*}
		&\underset{k\in\mathbb{N}}{\sup}\mathbb{E} \underset{t\in[0,T]}{\sup}\left \| u^{n}(t;\varphi_{k}) \right \|^{2}_{U_{1}} +\underset{k\in\mathbb{N}}{\sup}\mathbb{E} \underset{t\in[0,T]}{\sup}\left \| u^{n}_{t}(\varphi_{k}) \right \|^{2}_{h}+\underset{k\in\mathbb{N}}{\sup}\mathbb{E}\int _{0}^{T}[\left \| (A(s, u^{n}(t;\varphi_{k})) \right \|_{B^{*}}^{\frac{p}{p-1}}\nonumber
 \\&+\left \| f^{n}(s,u^{n}_{s}(\varphi_{k})) \right \|^{2}_{U_{1}}
+\left \| g^{n}(s,u^{n}_{t}(\varphi_{k})) \right \| ^{2}_{\mathscr{L}(U_{2},U_{1})}+\left \| u^{n}(t;\varphi_{k}) \right \|_{B}^{p} ]\text{d}s< \infty, \nonumber
	\end{align*}
and
\begin{align*}
		&\mathbb{E} \underset{t\in[0,T]}{\sup}\left \| u^{n}(t;\varphi)  \right \|^{2}_{U_{1}} +\mathbb{E} \underset{z\in[0,T]}{\sup}\left \| u^{n}_{t}(\varphi) \right \|^{2}_{h}+\mathbb{E}\int _{0}^{T}[\left \| (A(s,u^{n}(t;\varphi) ) \right \|_{B^{*}}^{\frac{p}{p-1}}\nonumber
 \\&+\left \| f^{n}(s,u^{n}_{t}(\varphi)) \right \|^{2}_{U_{1}}
+\left \| g^{n}(s,u^{n}_{t}(\varphi)) \right \| ^{2}_{\mathscr{L}(U_{2},U_{1})}+\left \|  u^{n}(t;\varphi) \right \|_{B}^{p} ]\text{d}s<  \infty. \nonumber
	\end{align*}
For $M > 0$,   we define the stopping time
\begin{align*}
    \tau^M_k &:= \inf \bigg\{ t \in [0, T ]; \|u^{n}(t;\varphi_{k})\|_{U_{1}} > M \bigg\} \wedge
    \inf \bigg\{ t \in [0, T ] : \int_0^t \|u^{n}(t;\varphi_{k})\|_B^p ds > M \bigg\} \\
    &~~~\wedge \inf \bigg\{ t \in [0, T ]: \|u^{n}(t;\varphi)\|_{U_{1}} > M \bigg\} \wedge
    \inf \bigg\{ t \in [0, T ] : \int_0^t \|u^{n}(t;\varphi)\|_B^p ds > M \bigg\}
    \wedge T,
\end{align*}
with the convention $\inf \emptyset = +\infty$. Then by the above apriori estimates for the solutions, we have
\begin{align}\label{iiiii1}
    \lim_{M\to\infty} \sup_{k\in\mathbb{N}} \mathbb{P}(\tau^M_k  < T ) = 0.
\end{align}
By It\^o's formula, Burkholder-Davis-Gundy inequality and Gronwall's inequality, we obtain
\begin{align}\label{iiiii2}
    \mathbb{E} \left[ \sup_{t\in[0,T]} \|u^{n}(t\wedge\tau^M_k;\varphi_{k}) - u^{n}(t\wedge\tau^M_k;\varphi)\|^{2}_{U_{1}} \right]
    \leq L_{M,T}\|\varphi_{k} - \varphi\|^{2}_{h},
\end{align}
which implies that for any $\epsilon > 0$ and $t\in [0,T]$, there exists a constant $L_M > 0$ such that
\begin{align}\label{iiiii7}
    &\mathbb{P}\left(\|u^{n}_{t}(\varphi_{k}) - u^{n}_{t}(\varphi)\|_{U_{1}} \geq \epsilon\right)\nonumber
    \\&\leq \mathbb{P}\left(\|u^{n}_{t}(\varphi_{k}) - u^{n}_{t}(\varphi)\|_{U_{1}} \geq \epsilon, \tau^M_k \geq T\right) + \mathbb{P}(\tau^M_k < T) \nonumber\\
    &\leq \frac{L_{M,T}}{\epsilon^2 } \mathbb{E} \left[ \|u^{n}_{t\wedge\tau^M_k}(\varphi_{k}) - u^{n}_{t\wedge\tau^M_k}(\varphi)\|^{2}_{h}\right] + \mathbb{P}(\tau^M_k < T) \\
    &\leq \frac{L_{M,T}\|\varphi_{k} - \varphi\|^{2}_{h}}{\epsilon^2 } + \sup_{k\in\mathbb{N}} \mathbb{P}(\tau^M_k < T).\nonumber
\end{align}
Applying\eqref{iiiii1}, we first let $k \to \infty$, followed by $M \to \infty$, to obtain that
\begin{align*}
    \lim_{k\to\infty} \|u^{n}_{t}(\varphi_{k}) - u^{n}_{t}(\varphi)\|_{h} = 0
\end{align*}
in probability for any $t \in [0, T]$. Furthermore, by $\underset{k\in\mathbb{N}}{\sup}\mathbb{E} \underset{t\in[0,T]}{\sup}\left \| u^{n}_{t}(\varphi_{k}) \right \|^{2}_{h}<\infty,$ it follows that
\begin{align}\label{iiiii8}
    \lim_{k\to\infty} \mathbb{E} \int_0^T \|u^{n}_{t}(\varphi_{k}) - u^{n}_{t}(\varphi)\|^{2}_{U_{1}} dt = 0.
\end{align}
Therefore,
\begin{align}\label{iiiii9}
    \lim_{k\to\infty} \mathbb{E} \int_0^T \|g(t, u^{n}_{t}(\varphi_{k})) - g(t, u^{n}_{t}(\varphi))\|_{\mathscr{L}(U_{2},U_{1})}^2 dt = 0.
\end{align}
By It\^o's formula, Burkholder-Davis-Gundy inequality and Young’s inequality, there exist a constant $a>0$ such that
\begin{align*}
    &\mathbb{E} \sup_{t\leq T \wedge \tau^M_k} e^{-at} \left \| u^{n}(t;\varphi_{k})- u^{n}(t;\varphi)\right \|^{2}_{U_{1}}
    \\&\leq L_{M}\|\varphi_{k} - \varphi\|^{2}_{h} + 2\mathbb{E} \sup_{t\leq T \wedge \sigma^M_n} \left| \int_0^t \langle e^{-as} u^{n}(s;\varphi_{k})- u^{n}(s;\varphi),  g(t, u^{n}_{t}(\varphi_{k})) - g(t, u^{n}_{t}(\varphi)) dW(s)\rangle_{U_{1}} \right| \\
    &\leq L_{M}\|\varphi_{k} - \varphi\|^{2}_{h} +\frac{1}{2}\mathbb{E} \sup_{t\leq T \wedge \tau^M_k} e^{-at} \left \| u^{n}(t;\varphi_{k})- u^{n}(t;\varphi)\right \|^{2}_{U_{1}}
     \\&~~~+L \mathbb{E} \int_0^{T \wedge \sigma^M_n} e^{-at} \|g(t, u^{n}_{t}(\varphi_{k})) - g(t, u^{n}_{t}(\varphi))\|_{\mathscr{L}(U_{2},U_{1})}^2 dt.
\end{align*}
From \eqref{iiiii8} and \eqref{iiiii9}, we conclude
\begin{align*}
    \lim_{k\to\infty} \mathbb{E} \sup_{t\leq T \wedge \tau^M_k} e^{-at} \left \| u^{n}(t;\varphi_{k})- u^{n}(t;\varphi)\right \|^{2}_{U_{1}} = 0.
\end{align*}
Arguing as in \eqref{iiiii7} again yields
\begin{align*}
    \lim_{k\to\infty}  \sup_{t\in[0, T]  } \left \| u^{n}(t;\varphi_{k})- u^{n}(t;\varphi)\right \|_{U_{1}} = 0
\end{align*}
in probability. Since $\underset{k\in\mathbb{N}}{\sup}\mathbb{E} \underset{t\in[0,T]}{\sup}\left \| u^{n}(t;\varphi_{k}) \right \|^{2}_{h}<\infty,$ we deduce that
\begin{align*}
    \lim_{k\to\infty} \mathbb{E} \sup_{t\in[0,T]} \|u^{n}(t;\varphi_{k})- u^{n}(t;\varphi)\|_{U_{1}}^{2} = 0,
\end{align*}
which implies
\begin{align*}
    \lim_{k\to\infty} |\mathbb{E}\mathcal{L}(u^{n}(t;\varphi_{k})|_{[0,T]})-\mathbb{E}\mathcal{L}( u^{n}(t;\varphi)|_{[0,T]})| = 0,
\end{align*}
for all $n\in \mathbb{N}$, i.e., the aforementioned convergence property holds uniformly in $n$. In addition, from the proof of statement 2), it follows that:
\begin{align*}
    \lim_{n\to\infty} |\mathbb{E}\mathcal{L}(u^{n}(t;\varphi_{k})|_{[0,T]})-\mathbb{E}\mathcal{L}( u(t;\varphi_{k})|_{[0,T]})| = 0,
\end{align*}
and \begin{align*}
    \lim_{n\to\infty} |\mathbb{E}\mathcal{L}(u^{n}(t;\varphi)|_{[0,T]})-\mathbb{E}\mathcal{L}( u(t;\varphi)|_{[0,T]})| = 0.
\end{align*}
Letting $n\to \infty$, we conclude that
\begin{align*}
    \lim_{k\to\infty} |\mathbb{E}\mathcal{L}(u(t;\varphi_{k})|_{[0,T]})-\mathbb{E}\mathcal{L}( u(t;\varphi)|_{[0,T]})| = 0.
\end{align*}
This completes the proof.

\end{enumerate}\quad $\Box$
\\~
\\\textbf{Remark 3.5.} Theorems 3.1 and 3.3 consider  two distinct classes of stochastic functional differential equations: one involving infinite delay in finite dimensions, and the other in infinite dimensions spaces. Our investigation centers on proving the existence and weak uniqueness of solutions under lower regularity assumptions, such as H$\ddot{\text{o}}$lder continuity of the coefficients. Moreover, for both finite delay and the case without delay, we show that analogous results can be obtained by setting  $h=0$ and appropriately adjusting the proof. This adjustment leads to a more concise derivation of the results.
~\\

Let $\tau \in (0,+\infty )$ be the delay, and define $\mathcal{H}   := C([-\tau , 0]; U_{1})$ with the norm $\left \| \varphi  \right \| _{h  } =\underset{-\tau\le\theta \le 0}{\sup}\left \|\varphi (\theta) \right \|_{U_{1}}  $, where $\mathcal{H}$  is the space of all continuous functions from $[-\tau , 0]$ into $U_{1}$. We now consider the following SFPDEs with finite delay
\begin{align}\label{j3}
		\begin{cases}
 \text{d}u(t)=(A(t,u(t))+f(t,u_{t})\text{d}t+g(t,u_{t})\text{d}W(t), \\
u_{0} =\varphi\in \mathcal{H},
\end{cases}
	\end{align}
where $u_{t} =u_{t} (\theta)=:\{u(t+\theta ), \theta\in[ -\tau,   0]\}$ and $f:\mathbb{R}^{+}\times\mathcal{H}  \to U_{1}$, $g:\mathbb{R}^{+}\times\mathcal{H}  \to \mathscr{L}(U_{2},U_{1})$ are two continuous mappings.
~\\
\\\textbf{Theorem 3.6.} \emph{Suppose that the assumptions \textbf{(H1)}$-$\textbf{(H5)} hold,  then the following statements hold:
\begin{enumerate}[1)]
        \item  For any initial value $\varphi \in \mathcal{H} $, system \eqref{j3} has  a probabilistic weak solution  $u(t)$ and a segment process $u_{t}$ to \eqref{j3} with $u_{0}=\varphi$;
        \item The weak solution $u(t)$  is unique in law; that is, any two such solutions with the same initial value have the same law;
		\item  Let $\{\varphi_{k}\}^{\infty}_{k=1}\in \mathcal{H}$ and $\varphi\in \mathcal{H}$ be a sequence with $\|\varphi_{k}-\varphi\|_{h}\to 0$. Let $u(t;\varphi)$ be the probabilistic weak solution of \eqref{j3} with the initial value $\varphi$. Then we have
\[\lim_{k\to\infty} |\mathbb{E}\mathcal{L}(u(t;\varphi_{k})|_{[0,T]})-\mathbb{E}\mathcal{L}( u(t;\varphi)|_{[0,T]})| = 0.\]
	\end{enumerate}}
For the proof of Theorem 3.6, the critical steps involve modifying part one from \eqref{z1} to
\begin{align*}
\left \| x^{n}_{s} \right \| _{h}&=\underset{-\tau < \Re \le 0}{\sup}\left |  x^{n} (s+\Re ) \right | \nonumber \\&\le\underset{-\tau < r \le s}{\sup}\left |  x^{n} (r) \right |
\\&\le\underset{-\tau < r \le 0}{\sup}\left |  x^{n} (r) \right |+\underset{0 \le r \le s}{\sup}\left |  x^{n} (r) \right | \nonumber
\\&\le\left \| \varphi  \right \|  _{h}+\underset{r\in[0,s]}{\sup}\left | x^{n}(r) \right |. \nonumber
\end{align*}
Let $\mathcal{C}\subseteq\mathcal{H}$ be a compact set ensuring that  $\mathbb{P}(\tau_{1} \leq T)$ becomes arbitrarily small. Define $\mathcal{C}$
  as the set of $\phi\in \mathcal{C}$ such that there exists an $r\in [0,\tau]$,
   \[\left\{\begin{matrix}
 \phi (\theta )=\varphi (\theta ), & \theta \in (-\tau ,-\tau+r],\\
 \phi (\theta )\in\mathbb{U},  & \theta \in [-\tau+r ,0],
\end{matrix}\right.\]
where $\varphi\in \mathcal{C}_{0}$ and $\mathcal{C}_{0}\subseteq\mathcal{H}$ is a compact set.

The remainder of the proof follows a procedure similar to those of Theorems 3.1 and 3.3; hence, we omit the specific details.
~\\

For the following systems without time delay:
\begin{align}\label{j5}
		\begin{cases}
 \text{d}u(t)=(A(t,u(t))+f(t,u(t))\text{d}t+g(t,u(t))\text{d}W(t), \\
u(0) =u_{0}\in U_{1},
\end{cases}
	\end{align}
it is necessary to adjust \textbf{(H2)} and \textbf{(H4)} as described below:
~\\
\\\textbf{(H2')} (Growth) For $A$, there exist constant $\alpha_{1}, M>0$ such that  for all $t\in\mathbb{R} ^{+}$  and $u\in B$, the following inequality holds:
 \begin{align*}
		\left \| A(t,u)\right \|_{B^{*}}^{\frac{p}{p-1}}   \le \alpha _{1}\left \| u\right \|_{B}^{p}+ M.
	\end{align*}
For the continuous functions $f$, $g$, there exist constant $\alpha _{1}$ and $M$ such that for all $t\in\mathbb{R} ^{+}$ and $u\in U_{1} $,
\begin{align*}
		\left \| f(t,u) \right \|_{U_{1}} \vee \left \| g(t,u ) \right \| _{\mathscr{L}(U_{2},U_{1})}  \le\alpha _{1}\left \| u\right \|_{U_{1}  }+M.
	\end{align*}
\textbf{(H4')}  There exist constants $\beta \in(0,1] $ and $\gamma  \in(\frac{1}{2},1 ]$. The map $A$ satisfies, for all $t\in\mathbb{R} ^{+}$ and $u,v\in B$,
\begin{align*}
		2_{B ^{\ast}} \langle A(t,u)-A(t,v), u-v\rangle _{B} \le \alpha _{1}\left \| u-v \right \| _{U_{1}}^{2},
	\end{align*}
and   the functions $f$, $g$ satisfy,  for all $t\in\mathbb{R} ^{+}$ and $u,v\in U_{1}  $ with $\left \| u\right \| _{U_{1}}\vee\left \| v \right \| _{U_{1}}\le M$,
\begin{align*}
\left \langle f(t,u)-f(t,v),u-v \right \rangle _{U_{1}}\le L_{M}\left \| u-v  \right \|^{\beta+1  } _{U_{1}},
\end{align*}
and
\begin{align*}
		\left \| g(t,u )- g(t,v ) \right \| _{\mathscr{L}(U_{2},U_{1})} \le L_{M}\left \|  u-v  \right \|^{\gamma  } _{U_{1}}.
	\end{align*}
\textbf{Theorem 3.7.} \emph{Consider \eqref{j5}. Suppose that the assumptions \textbf{(H1)}, \textbf{(H2')}, \textbf{(H3)}, \textbf{(H4')} and \textbf{(H5)} hold, then the following statements hold:
\begin{enumerate}[1)]
        \item  For any initial value $u_{0} \in U_{1} $, system \eqref{j5} has  a probabilistic weak solution  $u(t)$ and the weak solution $u(t)$  is unique in law, i.e., any two such solutions with the same initial value have the same law;
		\item  Let $\{u^{k}_{0}\}^{\infty}_{k=1}\in U_{1}$ and $u_{0}\in U_{1}$ be a sequence with $\|u^{k}_{0}-u_{0}\|_{U_{1}}\to 0$. Let $u(t;u_{0})$ be the probabilistic weak  solution of \eqref{j5} with the initial value $u_{0}$. Then we have
\[\lim_{k\to\infty} |\mathbb{E}\mathcal{L}(u(t;u^{k}_{0})|_{[0,T]})-\mathbb{E}\mathcal{L}( u(t;u_{0})|_{[0,T]})| = 0.\]
	\end{enumerate}}

Similarly, the proof of Theorem 3.7 follows a procedure similar to those of Theorems 3.1 and 3.3; hence, we omit the details.
\\~
\\\textbf{Remark 3.8.} In comparison with Theorem 3.1, the proof of Theorem 3.7 differs in two main aspects. First, regarding the construction of Lipschitz continuous approximating sequences ${f_n},{g_n}({F_n},{G_n})$: Since $U_{1}$ is a separable Hilbert space, we apply Proposition 3.3 of \cite{ref777} to construct  a sequence of bounded Lipschitz functions $\{f_n\},\{g_n\}(\{F_n\},\{G_n\})$ that converges uniformly to $f,g(F,G)$ on any compact subset of $U_{1}$. Second, in the definition of stopping times $\tau_{1}$: Let $\tau_{1}=\inf \{t\ge 0;u(t)\notin \mathbb{U}\}$, where $\mathbb{U}$ denotes a compact subset of $U_{1}$, so that the final conclusion follows from \eqref{r15r}.

\section{\textup{The weak averaging principle of  SFPDEs  with H$\ddot{\text{o}}$lder coefficients}}
Consider the following SFPDEs with infinite delay
\begin{align}\label{a1}
		\begin{cases}
 \text{d}u^{\varepsilon }(t)=(A(u^{\varepsilon }(t))+f(\frac{t}{\varepsilon } ,u^{\varepsilon }_{t})\text{d}t+g(\frac{t}{\varepsilon } ,u^{\varepsilon }_{t})\text{d}W(t), \\
u_{0} =\varphi^{\varepsilon }\in C^{h}_{U_{1}},
\end{cases}
	\end{align}
where $\varepsilon \in (0,1]$. Let $f^{\varepsilon}(t,\varphi )=f(\frac{t}{\varepsilon } ,\varphi)$ and $g^{\varepsilon}(t,\varphi )=g(\frac{t}{\varepsilon } ,\varphi)$, then \eqref{a1} can be written as
\begin{align}\label{a3}
		\begin{cases}
 \text{d}u^{\varepsilon }(t)=(A(u^{\varepsilon }(t))+f^{\varepsilon}(t, u^{\varepsilon }_{t})\text{d}t+g^{\varepsilon}(t,u^{\varepsilon }_{t})\text{d}W(t), \\
u_{0} =\varphi^{\varepsilon }\in C^{h}_{U_{1}}.
\end{cases}
	\end{align}
In order to conduct a more thorough analysis of the averaging principle applied to systems \eqref{a3} with H$\ddot{\text{o}}$lder coefficients, we introduce the following assumption:
~\\
\\\textbf{(H6)} There exist functions $\Phi _{1} $, $\Phi _{2} $ and $f^{*}\in C(C^{h}_{U_{1}},U_{1})$, $g^{*}\in C(C^{h}_{U_{1}},\mathscr{L}(U_{2},U_{1}) )$ such that for any $t\ge0$, $r>0$ and $\varphi \in C^{h}_{U_{1}}$,
\begin{align*}
		\frac{1}{r} \left \| \int_{t}^{t+r}[f(s,\varphi ) -f^{*}(\varphi )] \text{d}s\right \| _{U_{1}}\le \Phi _{1}(r)(\left \| \varphi  \right \|_{h}+M ),
	\end{align*}
\begin{align*}
		\frac{1}{r}\int_{t}^{t+r} \left \| g(s,\varphi ) -g^{*}(\varphi )\right \|^{2} _{\mathscr{L}(U_{2},U_{1}) } \text{d}s\le \Phi_{2} (r)(\left \| \varphi  \right \|^{2}_{h}+M ),
	\end{align*}
where $\Phi _{i} $ is decreasing, positive bounded functions and $\lim_{r \to \infty} \Phi _{i}(r)=0$ for $i=1,2$. In addition, we require that $g^{*}$ also satisfies condition \textbf{(H5)}.

Now we consider the following averaged equation:
\begin{align}\label{a2}
		\begin{cases}
 \text{d}u(t)=(A(u(t))+f^{*}(u_{t})\text{d}t+g^{*}(u_{t})\text{d}W(t), \\
u_{0} =\varphi^{* }\in C^{h}_{U_{1}}.
\end{cases}
	\end{align}
\textbf{Theorem 4.1.} \emph{Consider \eqref{a2}. If $A$, $f$ and $g$ satisfy \textbf{(H1)}$-$\textbf{(H6)}, the following statement holds: for any $\varphi^{*}\in C^{h}_{U_{1}} $, there exist a  probabilistic weak solution $u^{*}(t)$ and a segment process $u^{*}_{t}$ to \eqref{a2} with $u^{*}_{0}=\varphi^{*}$. Furthermore,  the weak solution $u^{*}(t)$  is unique in law.
\\\textbf{proof:}} According to  Theorem 3.2,  it suffices to verify that assumptions \textbf{(H3)}$-$\textbf{(H4)} hold for the coefficients $f^{*}$, $g^{*}$. Under the assumption \textbf{(H6)} we have
\begin{align*}
		\left \| f^{*}(\varphi) \right \|_{U_{1}} =  \left \|\lim_{r \to \infty}\frac{1}{r} \int_{t}^{t+r}f(s,\varphi ) \text{d}s\right \| _{U_{1}}\le \beta _{1}\left \| \varphi\right \|_{\mathcal{R}  }+M,
	\end{align*}
and
\begin{align*}
		\left \| f^{*}(\varphi)- f^{*}(\phi ) \right \|_{U_{1}}  &\le \left \| f^{*}(\varphi)- \frac{1}{r} \int_{0}^{r}f(s,\varphi ) \text{d}s \right \|_{U_{1}}+\left \| f^{*}(\phi )- \frac{1}{r} \int_{0}^{r}f(s,\phi ) \text{d}s \right \|_{U_{1}}
\\&~~~+\left \| \frac{1}{r} \int_{0}^{r}f(s,\varphi ) \text{d}s- \frac{1}{r} \int_{0}^{r}f(s,\phi ) \text{d}s \right \|_{U_{1}}
\\&\le\Phi _{1}(r)(\left \| \varphi  \right \|_{h}+M )+\Phi _{1}(r)(\left \| \phi  \right \|_{h}+M )+ L_{M}\left \| \varphi -\phi  \right \|^{\gamma  } _{h},
	\end{align*}
 As $r\to \infty $, we obtain for all $\varphi,\phi\in C^{h}_{U_{1}}  $ with $\left \| \varphi\right \| _{h}\vee\left \| \phi  \right \| _{h}\le M$,
 \begin{align*}
		\left \| f^{*}(\varphi)- f^{*}(\phi ) \right \|_{U_{1}}  \le L_{M}\left \| \varphi -\phi  \right \|^{\gamma  } _{h}.
	\end{align*}
The case of $g^{*}$ is similar.

For the assumption \textbf{(H4)}, we have for any $\varphi,\phi\in C^{h}_{U_{1}}$,
\begin{align*}
		&\left \langle f^{*}(\varphi)-f^{*}(\phi),\varphi(0)-\phi(0) \right \rangle _{U_{1}}
\\&\le \left \langle f^{*}(\varphi)-\frac{1}{r} \int_{0}^{r}f(s,\varphi ) \text{d}s,\varphi(0)-\phi(0) \right \rangle _{U_{1}}-\left \langle f^{*}(\phi)-\frac{1}{r} \int_{0}^{r}f(s,\phi) \text{d}s,\varphi(0)-\phi(0) \right \rangle _{U_{1}}
\\&~~~+\left \langle \frac{1}{r} \int_{0}^{r}f(s,\varphi ) \text{d}s-\frac{1}{r} \int_{0}^{r}f(s,\phi) \text{d}s,\varphi(0)-\phi(0) \right \rangle _{U_{1}}
\\&\le\Phi _{1}(r)(\left \| \varphi  \right \|_{h}+\left \| \phi  \right \|_{h}+M )\left \| \varphi -\phi  \right \| _{h}
+L_{M}\left \| \varphi -\phi  \right \|_{h}^{\gamma+1  } ,
	\end{align*}
and
\begin{align*}
	 &\left \| g^{*}(\varphi)-g^{*}(\phi) \right \|^{2} _{\mathscr{L}(U_{2},U_{1})}
\\&=\frac{1}{r} \int_{0}^{r}\left \| g^{*}(\varphi)-g^{*}(\phi) \right \|^{2} _{\mathscr{L}(U_{2},U_{1})}\text{d}s
\\&\le \frac{1}{r} \int_{0}^{r}[\left \| g^{*}(\varphi)-g(s,\varphi) \right \|^{2} _{\mathscr{L}(U_{2},U_{1})}+\left \| g^{*}(\phi)-g(s,\phi) \right \|^{2} _{\mathscr{L}(U_{2},U_{1})}
\\&~~~+\left \| g(s,\varphi)-g(s,\phi) \right \|^{2} _{\mathscr{L}(U_{2},U_{1})}]\text{d}s
\\&\le\Phi _{1}(r)(\left \| \varphi  \right \|_{h}+\left \| \phi  \right \|_{h}+M )+L_{M} \left \| \varphi -\phi  \right \|_{h}^{2\gamma  } .
	\end{align*}
Similarly, as $r\to \infty $, the coefficients $f$ and $g$ satisfy assumption \textbf{(H4)}.  This completes the proof. \quad $\Box $

To obtain the weak average principle of the above system, we continue to use the generalized coupling method and construct the auxiliary system as follows:
\begin{align*}
    \begin{cases}
        \text{d}u^{\varepsilon,n,\zeta}(t) =  \left(A(u^{\varepsilon,n,\zeta}(t)) + f^{\varepsilon}_{n}(t,u^{\varepsilon,n,\zeta}_{t})+ \zeta(u^{\varepsilon}(t) - u^{\varepsilon,n,\zeta}(t))\mathbf{1} _{t\le\tau}\right) \text{d}t
         + g^{\varepsilon}_{n}(t,u^{\varepsilon,n,\zeta}_{t}) \text{d}W(t)
       \text{d}t, \\
        u_{0} =\varphi^{\varepsilon}\in C^{h}_{U_{1}},
    \end{cases}
\end{align*}
and
\begin{align*}
    \begin{cases}
        \text{d}u^{*,n,\zeta}(t) =  \left(A(u^{*,n,\zeta}(t)) + f^{*}_{n}(u^{*,n,\zeta}_{t})+ \zeta(u^{*}(t) - u^{*,n,\zeta}(t))\mathbf{1} _{t\le\tau}\right) \text{d}t
         + g^{*}_{n}(u^{*,n,\zeta}_{t}) \text{d}W(t)
       \text{d}t, \\
        u_{0} =\varphi^{*}\in C^{h}_{U_{1}},
    \end{cases}
\end{align*}
where families $\{f^{\varepsilon}_{n}\}, \{g^{\varepsilon}_{n}\}$($\{f^{*}_{n}\}, \{g^{*}_{n}\}$) satisfy
\begin{enumerate}[(\textbf{I})]
		\item $f^{\varepsilon}_{n}\to f^{\varepsilon}, g^{\varepsilon}_{n}\to f^{\varepsilon}$($f^{*}_{n}\to f^{*}, g^{*}_{n}\to f^{*}$) as $n\to \infty $  uniformly on each compact subset of $ C^{h}_{U_{1}}$;
\end{enumerate}
\begin{enumerate}[(\textbf{II})]
        \item $\{f^{\varepsilon}_{n}\}, \{g^{\varepsilon}_{n}\}$($\{f^{*}_{n}\}, \{g^{*}_{n}\}$) satisfy the conditions \textbf{(H1)}-\textbf{(H5)}, that is,  the constants are independent of $n$.
	\end{enumerate}
\begin{enumerate}[(\textbf{III})]
        \item the functions $\{f^{\varepsilon}_{n}\}, \{g^{\varepsilon}_{n}\}$($\{f^{*}_{n}\}, \{g^{*}_{n}\}$) are Lipschitz continuous on each bounded subset of $C^{h}_{U_{1}}$.
	\end{enumerate}
The above systems should be understood as a controlled version of
\begin{align}
    \begin{cases}
        \text{d}u^{\varepsilon,n}(t) = \left((A(u^{\varepsilon,n}(t)) + f^{\varepsilon}_{n}(t,u^{\varepsilon,n}_{t})\right) \text{d}t + g^{\varepsilon}_{n}(t,u^{\varepsilon,n}_{t}) \text{d}W(t), \\
        u_{0} =\varphi^{\varepsilon}\in C^{h}_{U_{1}},
    \end{cases} \label{e1e}
\end{align}
and
\begin{align}
    \begin{cases}
        \text{d}u^{*,n}(t) = \left((A(u^{*,n}(t)) + f^{*}_{n}(u^{*,n}_{t})\right) \text{d}t + g^{*}_{n}(u^{*,n}_{t}) \text{d}W(t), \\
        u_{0} =\varphi^{*}\in C^{h}_{U_{1}}.
    \end{cases} \label{e2e}
\end{align}

In the following, we first study the related properties of solution $u^{\varepsilon,n}(t)$ of system  \eqref{e1e} and solution $u^{*,n}(t)$ of system  \eqref{e2e}. It is worth noting that the coefficients of system \eqref{e1e} and \eqref{e2e} satisfy the Lipschitz continuity condition, so $u^{\varepsilon,n}(t), u^{*,n}(t)$  are the unique strong solutions of these systems. For any given function $\Psi $, define a piecewise function $\hat{\Psi } $  as follows:
\begin{align}\label{a5}
	 \hat{\Psi  } (t)=\left\{\begin{matrix}
 \Psi(t)  & t< 0,\\
  \Psi (0)&t\in [0,d),\\
  ...&...\\
\Psi (kd)&t\in [kd,(k+1)d),\\
   ...&...
\end{matrix}\right.
	\end{align}
 where $k\in\mathbb{N} ^{+}$ and $d$ be a fixed constant. To further investigate the long-term asymptotic behavior of systems \eqref{a3} and \eqref{a2}, it is necessary to introduce the following lemma:
\\~
\\\textbf{Lemma 4.2.} \emph{If $A$, $f$ and $g$ satisfy \textbf{(H1)}$-$\textbf{(H6)}, the following statements hold:  for any $T>0$ and  $\varphi ^{\varepsilon }, \varphi ^{* }\in C^{h}_{U_{1}}$,
\begin{enumerate}[(1)]
		\item
\begin{align*}
		\mathbb{E} \int_{0}^{T} \left \| u^{\varepsilon,n }(s;\varphi ^{\varepsilon }) -\hat{u}^{\varepsilon,n }(s;\varphi ^{\varepsilon }) \right \|^{2}_{U_{1}} \text{d}s\le L_{T}(\left \| \varphi ^{\varepsilon } \right \|^{2}_{h}+1 )d^{\frac{1}{2} } ,
	\end{align*}
and
\begin{align*}
		\mathbb{E} \int_{0}^{T} \left \| u^{*,n }(s;\varphi ^{* }) -\hat{u}^{*,n}(s;\varphi ^{*}) \right \|^{2}_{U_{1}} \text{d}s\le L_{T}(\left \| \varphi ^{* } \right \|^{2}_{h}+1 )d^{\frac{1}{2} };
	\end{align*}
	\item \begin{align*}
		\mathbb{E} \int_{0}^{T} \left \| u^{\varepsilon,n }_{s}(\varphi ^{\varepsilon }) -\hat{u}^{\varepsilon,n }_{s}(\varphi ^{\varepsilon }) \right \|^{2}_{h} \text{d}s\le L_{T}(\left \| \varphi ^{\varepsilon } \right \|^{2}_{h}+1 )d^{\frac{1}{2} } ,
	\end{align*}
and
\begin{align*}
		\mathbb{E} \int_{0}^{T} \left \| u^{*,n }_{s}(\varphi ^{* }) -\hat{u}^{*,n}_{s}(\varphi ^{*}) \right \|^{2}_{h} \text{d}s\le L_{T}(\left \| \varphi ^{* } \right \|^{2}_{h}+1 )d^{\frac{1}{2} },
	\end{align*}
	\end{enumerate}
where  $u^{\varepsilon,n }(s;\varphi ^{\varepsilon })$ $(u^{\varepsilon,n }_{s}(\varphi ^{\varepsilon }))$ is the solution(the solution map) of \eqref{a3} with the initial value $u^{\varepsilon,n }_{0}=\varphi ^{\varepsilon }$ and
$u^{*,n }(s;\varphi ^{*})$ $(u^{*,n }_{s}(\varphi ^{*}))$ is the solution(the solution map) of \eqref{a2} with the initial value $u^{*,n }_{0}=\varphi ^{*}$.
~\\
\\\textbf{proof of (1):}}  Similar to the proof of Lemma 4.4 in \cite{ref2}, the above results can be obtained using \textbf{(H1)}$-$\textbf{(H6)}, \eqref{p7}, It$\hat{\text{o}} $ formula, Burkholder-Davis-Gundy inequality and  Young's inequality.
~\\
 \\\emph{\textbf{proof of (2):}} By \eqref{a5}, we obtain
 \begin{align}\label{a6}
		&\mathbb{E} \int_{0}^{T} \left \| u^{\varepsilon,n }_{s}(\varphi ^{\varepsilon }) -\hat{u}^{\varepsilon,n }_{s}(\varphi ^{\varepsilon }) \right \|^{2}_{h} \text{d}s \nonumber
\\&=\mathbb{E} \int_{0}^{T} \underset{\theta  \in(-\infty ,0]}{\sup}e^{2h\theta  } \left \| u^{\varepsilon,n }(s+\theta  ;\varphi ^{\varepsilon }) -\hat{u}^{\varepsilon,n }(s+\theta  ;\varphi ^{\varepsilon }) \right \|^{2}_{U_{1}} \text{d}s \nonumber
\\&\le\mathbb{E} \underset{\theta   \in(-\infty ,0]}{\sup}\int_{\theta  }^{T+\theta  } \left \| u^{\varepsilon,n }(z;\varphi ^{\varepsilon }) -\hat{u}^{\varepsilon,n }(z;\varphi ^{\varepsilon }) \right \|^{2}_{U_{1}} \text{d}z
\\&\le\mathbb{E} \underset{\theta   \in(-\infty ,0]}{\sup}\int_{\theta  }^{0} \left \| u^{\varepsilon,n }(z;\varphi ^{\varepsilon }) -\hat{u}^{\varepsilon,n }(z;\varphi ^{\varepsilon }) \right \|^{2}_{U_{1}} \text{d}z+\mathbb{E} \int_{0}^{T} \left \| u^{\varepsilon,n }(z;\varphi ^{\varepsilon }) -\hat{u}^{\varepsilon,n }(z;\varphi ^{\varepsilon }) \right \|^{2}_{U_{1}} \text{d}z \nonumber
\\&\le L_{T}(\left \| \varphi ^{\varepsilon } \right \|^{2}_{h}+1 )d^{\frac{1}{2} }. \nonumber
	\end{align}
It follows from the same steps as \eqref{a6} that
\begin{align*}
		\mathbb{E} \int_{0}^{T} \left \| u^{*,n }_{s}(\varphi ^{* }) -\hat{u}^{*,n}_{s}(\varphi ^{*}) \right \|^{2}_{h} \text{d}s\le L_{T}(\left \| \varphi ^{* } \right \|^{2}_{h}+1 )d^{\frac{1}{2} }.
	\end{align*}
$\Box$
~\\

Now we establish  the weak averaging principle for SFPDEs  with infinite delay and  H$\ddot{\text{o}}$lder coefficients.
\\~
\\\textbf{Theorem 4.3.} \emph{Consider \eqref{a3} and \eqref{a2} under the assumptions \textbf{(H1)}-\textbf{(H6)}. For any initial values $\varphi ^{\varepsilon }, \varphi ^{*}\in C^{h}_{U_{1}}$ and $T > 0$, suppose further that $\lim_{\varepsilon  \to 0} \mathbb{E} \left \|\varphi ^{\varepsilon } -\varphi ^{*} \right \|^{2}_{h}=0$. Then, we have
\begin{align}\label{c1}
		\lim_{\varepsilon  \to 0}|\mathbb{E}\mathcal{L}(u^{\varepsilon }(t;\varphi ^{\varepsilon })| _{[0,T]}) - \mathbb{E}\mathcal{L}(u^{* }(t;\varphi ^{*})|_{[0,T]})|= 0,
	\end{align}
where $\mathcal{L}$ is  an arbitrary real-valued, bounded continuous function defined on $C([0, T], U_{1})$.
\\\textbf{proof:}} Given any compact subset  $\mathcal{C} \subset C^{h}_{U_{1}}$,  denote
\[
\tau^{*}=\inf \{t\ge 0;u_{t}^{\varepsilon}\notin \mathcal{C}~ \text{and} ~ u_{t}^{*}\notin \mathcal{C}\},\quad
 \Re^{\varepsilon}_{f}=\sup_{t\ge0}\sup_{\varphi\in \mathcal{C}}\| f^{\varepsilon}_{n}(t,\varphi)-f^{\varepsilon}(t,\varphi)\|_{U_{1}},
 \]
 \[ \Re^{*}_{f}=\sup_{\varphi\in \mathcal{C}}\| f^{*}_{n}(\varphi)-f^{*}(\varphi)\|_{U_{1}},\quad
\Re^{\varepsilon}_{g}=\sup_{t\ge0}\sup_{\varphi\in \mathcal{C}}\| g^{\varepsilon}_{n}(t,\varphi)-g^{\varepsilon}(t,\varphi)\|_{\mathscr{L}(U_{2},U_{1})},
\]
\[\Re^{*}_{g}=\sup_{\varphi\in \mathcal{C}}\| g^{*}_{n}(\varphi)-g^{*}(\varphi)\|_{\mathscr{L}(U_{2},U_{1})},\quad\Re(n)=\max\{(\Re^{\varepsilon}_{f})^{\frac{1}{\beta}},(\Re^{*}_{f})^{\frac{1}{\beta}} ,(\Re^{\varepsilon}_{g})^{\frac{1}{\gamma}}, (\Re^{*}_{g})^{\frac{1}{\gamma}}  \}.
\]
Since $\mathcal{C}$ is compact, the sequence $(f^{\varepsilon}_{n},g^{\varepsilon}_{n},f^{*}_{n},g^{*}_{n})$ converges to $(f^{\varepsilon},g^{\varepsilon},f^{*},g^{*})$ for each $\varphi\in \mathcal{C}$, and $(f^{\varepsilon},g^{\varepsilon},f^{*},g^{*})$ is
uniformly continuous in $\varphi$. Therefore, we deduce that  $\lim_{n \to \infty}\Re(n)=0 $, since all the estimates are time-independent. That is, for any $\varepsilon>0$, there exists $N$ such that for all $n>N$, $\Re(n)<\varepsilon$.

Applying It$\hat{\text{o}} $  formula to $u^{\varepsilon,n }(t;\varphi ^{\varepsilon }) -u^{*,n }(t;\varphi ^{*})$, and using  \textbf{(H4)}, Young inequality and Burkholder-Davis-Gundy inequality, we obtain
 \begin{align*}
		&\mathbb{E}\underset{t\in [0,\tau^{*}\wedge T]}{\sup}\left \| u^{\varepsilon,n }(t;\varphi ^{\varepsilon }) -u^{*,n }(t;\varphi ^{*}) \right \|^{2}_{U_{1}}\nonumber
\\&=\left \| \varphi ^{\varepsilon }(0) -\varphi ^{*}(0) \right \|^{2}_{U_{1}} +\mathbb{E}\underset{t\in [0,\tau^{*}\wedge T]}{\sup}\{\int_{0}^{t}[2_{B ^{\ast}} \langle A(u^{\varepsilon,n }(s))- A(u^{*,n }(s)),u^{\varepsilon,n }(s)- u^{*,n }(s)\rangle _{B}\nonumber
\\&~~~+2\langle f^{\varepsilon}_{n}(s,u^{\varepsilon,n }_{s})- f^{*}_{n}(u^{*,n }_{s}),u^{\varepsilon,n }(s)- u^{*,n }(s)\rangle_{U_{1}}+\left \| g^{\varepsilon}_{n}(s,u^{\varepsilon,n }_{s})- g^{*}_{n}(u^{*,n }_{s})\right \|^{2}_{\mathscr{L}(U_{2},U_{1})}]\text{d}s
\\&~~~+2\int_{0}^{t} \left \langle u^{\varepsilon,n }(s)- u^{*,n }(s), [g^{\varepsilon}_{n}(s,u^{\varepsilon,n }_{s})- g^{*}_{n}(u^{*,n }_{s}) ]\text{d}W(s)\right \rangle_{U_{1}}\}  \nonumber
\\&\le \left \| \varphi ^{\varepsilon }(0) -\varphi ^{*}(0) \right \|^{2}_{U_{1}} +\alpha_{1}\mathbb{E}\int_{0}^{\tau^{*}\wedge T}[\left \| u^{\varepsilon,n }(s)- u^{*,n }(s) \right \| _{U_{1}}^{2 }\text{d}s  \nonumber
\\&~~~+\mathbb{E}\underset{t\in [0,\tau^{*}\wedge T]}{\sup}\int_{0}^{t}2\langle f^{\varepsilon}_{n}(s,u^{\varepsilon,n }_{s})- f^{*}_{n}(u^{*,n }_{s}),u^{\varepsilon,n }(s)- u^{*,n }(s)\rangle_{U_{1}}
\\&~~~+73\left \| g^{\varepsilon}_{n}(s,u^{\varepsilon,n }_{s})- g^{*}_{n}(u^{*,n }_{s})\right \|^{2}_{\mathscr{L}(U_{2},U_{1})}]\text{d}s
+\frac{1}{2} \mathbb{E}\underset{t\in [0,\tau^{*}\wedge T]}{\sup}\left \| u^{\varepsilon,n }(t;\varphi ^{\varepsilon }) -u^{*,n }(t;\varphi ^{*}) \right \|^{2}_{U_{1}},
	\end{align*}
which implies
 \begin{align}\label{a7}
		&\mathbb{E}\underset{t\in [0,\tau^{*}\wedge T]}{\sup}\left \| u^{\varepsilon,n }(t;\varphi ^{\varepsilon }) -u^{*,n }(t;\varphi ^{*}) \right \|^{2}_{U_{1}}\nonumber
\\&\le 2\left \| \varphi ^{\varepsilon } -\varphi ^{*} \right \|^{2}_{h} +2\alpha_{1}\mathbb{E}\int_{0}^{\tau^{*}\wedge T}[\left \| u^{\varepsilon,n }(s)- u^{*,n }(s) \right \| _{U_{1}}^{2 }\nonumber
\\&~~~+146\left \| g^{\varepsilon}_{n}(s,u^{\varepsilon,n }_{s})- g^{*}_{n}(u^{*,n }_{s})\right \|^{2}_{\mathscr{L}(U_{2},U_{1})}]\text{d}s
\\&~~~+4\mathbb{E}\underset{t\in [0,\tau^{*}\wedge T]}{\sup}\int_{0}^{t}\langle f^{\varepsilon}_{n}(s,u^{\varepsilon,n }_{s})- f^{*}_{n}(u^{*,n }_{s}),u^{\varepsilon,n }(s)- u^{*,n }(s)\rangle_{U_{1}}\text{d}s. \nonumber
	\end{align}
For the drift coefficient,
\begin{align}\label{a8}
	&\mathbb{E}\underset{t\in [0,\tau^{*}\wedge T]}{\sup}\int_{0}^{t}\langle f^{\varepsilon}_{n}(s,u^{\varepsilon,n }_{s})- f^{*}_{n}(u^{*,n }_{s}),u^{\varepsilon,n }(s)- u^{*,n }(s)\rangle_{U_{1}}\text{d}s\nonumber
\\&\le \mathbb{E}\underset{t\in [0,\tau^{*}\wedge T]}{\sup}\int_{0}^{t}\langle f^{\varepsilon}_{n}(s,u^{\varepsilon,n }_{s})- f^{\varepsilon}_{n}(s,u^{*,n }_{s}),u^{\varepsilon,n }(s)- u^{*,n }(s)\rangle_{U_{1}}\text{d}s
\\&~~~+\mathbb{E}\underset{t\in [0,\tau^{*}\wedge T]}{\sup}\int_{0}^{t}\langle f^{\varepsilon}_{n}(s,u^{*,n }_{s})- f^{*}_{n}(u^{*,n }_{s}),u^{\varepsilon,n }(s)- u^{*,n }(s)\rangle_{U_{1}}\text{d}s\nonumber
\\&:=\mathfrak{f} _{1}+\mathfrak{f}  _{2}. \nonumber
	\end{align}
By (\textbf{III}) and \eqref{z1}, we have
\begin{align}\label{a9}
	\mathfrak{f} _{1}&=\mathbb{E}\underset{t\in [0,\tau^{*}\wedge T]}{\sup}\int_{0}^{t}\langle f^{\varepsilon}_{n}(s,u^{\varepsilon,n }_{s})- f^{\varepsilon}_{n}(s,u^{*,n }_{s}),u^{\varepsilon,n }(s)- u^{*,n }(s)\rangle_{U_{1}}\text{d}s  \nonumber
\\&\le L_{M}\mathbb{E}\underset{t\in [0,\tau^{*}\wedge T]}{\sup}\int_{0}^{t}\left \| u^{\varepsilon,n }_{s} -u^{*,n }_{s} \right \|^{2 }_{h}\text{d}s
\\&\le L_{M}\int_{0}^{\tau^{*}\wedge T}\underset{z\in [0,s]}{\sup}\mathbb{E}\left \| u^{\varepsilon,n }(z) -u^{*,n }(z) \right \|^{2 }_{U_{1}}\text{d}s+L_{M,T}\left \| \varphi ^{\varepsilon } -\varphi ^{*} \right \| ^{2 }_{h}. \nonumber
	\end{align}
For $\mathfrak{f} _{2}$,
\begin{align}\label{a10}
	\mathfrak{f} _{2}&=\mathbb{E}\underset{t\in [0,\tau^{*}\wedge T]}{\sup}\int_{0}^{t}\langle f^{\varepsilon}_{n}(s,u^{*,n }_{s})- f^{*}_{n}(u^{*,n }_{s}),u^{\varepsilon,n }(s)- u^{*,n }(s)\rangle_{U_{1}}\text{d}s\nonumber
\\&\le \mathbb{E}\underset{t\in [0,\tau^{*}\wedge T]}{\sup}\int_{0}^{t}\langle f^{\varepsilon}_{n}(s,u^{*,n }_{s})- f^{*}_{n}(u^{*,n }_{s}),u^{\varepsilon,n }(s)- \hat{u} ^{\varepsilon,n}(s)\rangle_{U_{1}}\text{d}s\nonumber
\\&~~~+\mathbb{E}\underset{t\in [0,\tau^{*}\wedge T]}{\sup}\int_{0}^{t}\langle f^{\varepsilon}_{n}(s,u^{*,n }_{s})- f^{*}_{n}(u^{*,n }_{s}), \hat{u} ^{\varepsilon,n}(s)-\hat{u} ^{*,n}(s)\rangle_{U_{1}}\text{d}s
\\&~~~+\mathbb{E}\underset{t\in [0,\tau^{*}\wedge T]}{\sup}\int_{0}^{t}\langle f^{\varepsilon}_{n}(s,u^{*,n }_{s})- f^{*}_{n}(u^{*,n }_{s}), \hat{u} ^{*,n}(s)-u ^{*,n}(s)\rangle_{U_{1}}\text{d}s\nonumber
\\&:=\mathfrak{f}^{1}_{2}+\mathfrak{f} ^{2}_{2}+\mathfrak{f} ^{3}_{2}. \nonumber
	\end{align}
By \textbf{(H3)}, \eqref{p7}, \eqref{a5}, H\"older's inequality and Lemma 4.2, we get
\begin{align}\label{a11}
	\mathfrak{f}^{1} _{2}&=\mathbb{E}\underset{t\in [0,\tau^{*}\wedge T]}{\sup}\int_{0}^{t}\langle f^{\varepsilon}_{n}(s,u^{*,n }_{s})- f^{*}_{n}(u^{*,n }_{s}),u^{\varepsilon,n }(s)- \hat{u} ^{\varepsilon,n}(s)\rangle_{U_{1}}\text{d}s\nonumber
\\&\le\mathbb{E}\int_{0}^{\tau^{*}\wedge T}[\left \| f^{\varepsilon}_{n}(s,u^{*,n }_{s})\right \|_{U_{1}} +\left \| f^{*}_{n}(u^{*,n }_{s}) \right \|_{U_{1}}] \left \| u^{\varepsilon,n }(s)- \hat{u} ^{\varepsilon,n}(s) \right \|_{U_{1}} \text{d}s
\\&\le (\mathbb{E}\int_{0}^{\tau^{*}\wedge T}2L_{M}(\left \| u^{*,n }_{s}\right \|^{2}_{U_{1}} +M) \text{d}s)^{\frac{1}{2} }(\mathbb{E}\int_{0}^{\tau^{*}\wedge T} \left \| u^{\varepsilon,n }(s)- \hat{u} ^{\varepsilon,n}(s) \right \|^{2}_{U_{1}} \text{d}s)^{\frac{1}{2} }\nonumber
\\&\le L_{T,M}(\left \| \varphi ^{\varepsilon }  \right \|^{2}_{h} +\left \| \varphi ^{*} \right \|^{2}_{h} +1)d^{\frac{1}{4} }. \nonumber
	\end{align}
Similarly, for $\mathfrak{f}^{3} _{2}$,
\begin{align}\label{a12}
	\mathfrak{f}^{3} _{2}&=\mathbb{E}\underset{t\in [0,\tau^{*}\wedge T]}{\sup}\int_{0}^{t}\langle f^{\varepsilon}_{n}(s,u^{*,n }_{s})- f^{*}_{n}(u^{*,n }_{s}), \hat{u} ^{*}(s)-u ^{*}(s)\rangle_{U_{1}}\text{d}s
\\&\le L_{T,M}(\left \| \varphi ^{*} \right \|^{2}_{h} +1)d^{\frac{1}{4} }. \nonumber
	\end{align}
Next, the key problem is to estimate $\mathfrak{f}^{2} _{2}$:
\begin{align}\label{a13}
	\mathfrak{f} ^{2}_{2}&=\mathbb{E}\underset{t\in [0,\tau^{*}\wedge T]}{\sup}\int_{0}^{t}\langle f^{\varepsilon}_{n}(s,u^{*,n }_{s})- f^{*}_{n}(u^{*,n }_{s}), \hat{u} ^{\varepsilon,n}(s)-\hat{u} ^{*,n}(s)\rangle_{U_{1}}\text{d}s \nonumber
\\&\le\mathbb{E}\underset{t\in [0,\tau^{*}\wedge T]}{\sup}\int_{0}^{t}\langle f^{\varepsilon}_{n}(s,u^{*,n }_{s})- f^{\varepsilon}_{n}(s,\hat{u} ^{*,n}_{s}), \hat{u} ^{\varepsilon,n}(s)-\hat{u} ^{*,n}(s)\rangle_{U_{1}}\text{d}s \nonumber
\\&~~~+\mathbb{E}\underset{t\in [0,\tau^{*}\wedge T]}{\sup}\int_{0}^{t}\langle f^{\varepsilon}_{n}(s,\hat{u} ^{*,n}_{s})-f^{*}_{n}(\hat{u} ^{*,n}_{s}), \hat{u} ^{\varepsilon,n}(s)-\hat{u} ^{*,n}(s)\rangle_{U_{1}}\text{d}s
\\&~~~+\mathbb{E}\underset{t\in [0,\tau^{*}\wedge T]}{\sup}\int_{0}^{t}\langle f^{*}_{n}(\hat{u} ^{*,n}_{s})-f^{*}_{n}(u^{*,n }_{s}), \hat{u} ^{\varepsilon,n}(s)-\hat{u} ^{*,n}(s)\rangle_{U_{1}}\text{d}s \nonumber
\\&:=\mathfrak{f} ^{2,1}_{2}+\mathfrak{f} ^{2,2}_{2}+\mathfrak{f} ^{2,3}_{2}. \nonumber
	\end{align}
By (\textbf{III}), \eqref{p7}, \eqref{a5}, H\"older's inequality, Jensen's inequality, Theorem 4.1 and Lemma 4.2, we get
\begin{align}\label{a15}
	\mathfrak{f} ^{2,1}_{2}&=\mathbb{E}\underset{t\in [0,\tau^{*}\wedge T]}{\sup}\int_{0}^{t}\langle f^{\varepsilon}_{n}(s,u^{*,n }_{s})- f^{\varepsilon}_{n}(s,\hat{u} ^{*,n}_{s}), \hat{u} ^{\varepsilon,n}(s)-\hat{u} ^{*,n}(s)\rangle_{U_{1}}\text{d}s\nonumber
\\&\le \mathbb{E}\int_{0}^{\tau^{*}\wedge T}\left \| f^{\varepsilon}_{n}(s,u^{*,n }_{s})- f^{\varepsilon}_{n}(s,\hat{u} ^{*,n}_{s}) \right \| _{U_{1}}\left \| \hat{u} ^{\varepsilon,n}(s)-\hat{u} ^{*,n}(s) \right \|_{U_{1}}\text{d}s\nonumber
\\&\le L_{T}\mathbb{E}\int_{0}^{\tau^{*}\wedge T}\left \| u^{*,n }_{s}- \hat{u} ^{*,n}_{s}\right \| _{h}\left \| \hat{u} ^{\varepsilon,n}(s)-\hat{u} ^{*,n}(s) \right \|_{U_{1}}\text{d}s
\\&\le L_{T}(\mathbb{E}\int_{0}^{\tau^{*}\wedge T}\left \| u^{*,n }_{s}- \hat{u} ^{*,n}_{s} \right \|_{h}^{2} \text{d}s )^{\frac{1 }{2} }(\int_{0}^{\tau^{*}\wedge T}\mathbb{E}\left \|\hat{u} ^{\varepsilon,n}(s)-\hat{u} ^{*,n}(s) \right \|_{U_{1}}^{2  } \text{d}s )^{\frac{1 }{2} }\nonumber
\\&\le L_{T}(\left \| \varphi ^{*} \right \|_{h}^{2} +1)d^{\frac{1 }{4} }(\int_{0}^{\tau^{*}\wedge T}[\mathbb{E}\left \|\hat{u} ^{\varepsilon,n}(s) \right \|^{2}_{U_{1}}+\mathbb{E}\left \|\hat{u} ^{*,n}(s) \right \|^{2}_{U_{1}}] \text{d}s )^{\frac{1 }{2} }\nonumber
\\&\le L_{T}(\left \| \varphi ^{*} \right \|_{h}^{2}+\left \| \varphi ^{\varepsilon} \right \|_{h}^{2} +1)d^{\frac{1 }{4} }. \nonumber
	\end{align}
Similarly, for $\mathfrak{f}^{2,3} _{2}$,
\begin{align}\label{a16}
	\mathfrak{f} ^{2,3}_{2}&=\mathbb{E}\underset{t\in [0,\tau^{*}\wedge T]}{\sup}\int_{0}^{t}\langle f^{*}_{n}(\hat{u} ^{*,n}_{s})-f^{*}_{n}(u^{*,n }_{s}), \hat{u} ^{\varepsilon,n}(s)-\hat{u} ^{*,n}(s)\rangle_{U_{1}}\text{d}s \nonumber
\\&\le L_{T}(\left \| \varphi ^{*} \right \|_{h}^{2}+\left \| \varphi ^{\varepsilon} \right \|_{h}^{2} +1)d^{\frac{1 }{4} }.
	\end{align}
In the following step, we will use the time discretization technique to deal with $\mathfrak{f}^{2,2} _{2}$. Let $\left [ t \right ] $ denote the integer part of $t$, then note that
\begin{align}\label{a17}
	\mathfrak{f} ^{2,2}_{2}&=\mathbb{E}\underset{t\in [0,\tau^{*}\wedge T]}{\sup}\int_{0}^{t}\langle f^{\varepsilon}_{n}(s,\hat{u} ^{*,n}_{s})-f^{*}_{n}(\hat{u} ^{*,n}_{s}), \hat{u} ^{\varepsilon,n}(s)-\hat{u} ^{*,n}(s)\rangle_{U_{1}}\text{d}s  \nonumber
\\&=\mathbb{E}\underset{t\in [0,\tau^{*}\wedge T]}{\sup}\sum_{m=0}^{\left [ \frac{t}{d}  \right ]-1 } \int_{md}^{(m+1)d}\langle f^{\varepsilon}_{n}(s,u ^{*,n}_{md})-f^{*}_{n}(u ^{*,n}_{md}), u ^{\varepsilon,n}(md)-u ^{*,n}(md)\rangle_{U_{1}}\text{d}s  \nonumber
\\&~~~+\int_{\left [ \frac{t}{d}  \right ]d}^{t}\langle f^{\varepsilon}_{n}(s,u ^{*,n}_{\left [ t/d  \right ]d})-f^{*}_{n}(u ^{*,n}_{\left [ t/d \right ]d}), u ^{\varepsilon,n}(\left [ t/d  \right ]d)-u ^{*,n}(\left [ t/d  \right ]d)\rangle_{U_{1}}\text{d}s  \nonumber
\\&\le \mathbb{E}\underset{t\in [0,\tau^{*}\wedge T]}{\sup}\sum_{m=0}^{\left [ \frac{t}{d}  \right ]-1 } \left \|\int_{md}^{(m+1)d} f^{\varepsilon}_{n}(s,u ^{*,n}_{md})-f^{*}_{n}(u ^{*,n}_{md})  \text{d}s\right \|_{U_{1}} \left \| u ^{\varepsilon,n}(md)-u ^{*,n}(md) \right \| _{U_{1}}\nonumber
\\&~~~+L_{T}(\left \| \varphi ^{*} \right \|_{h}^{2}+\left \| \varphi ^{\varepsilon} \right \|_{h}^{2} +1)d
\\&\le\sum_{m=0}^{\left [ \frac{\tau^{*}\wedge T}{d}  \right ]-1 } (\left \|\mathbb{E}\int_{md}^{(m+1)d} f^{\varepsilon}_{n}(s,u ^{*,n}_{md})-f^{*}_{n}(u ^{*,n}_{md}) \text{d}s\right \|^{2}_{U_{1}})^{\frac{1}{2} } (\mathbb{E}\left \| u ^{\varepsilon,n}(md)-u ^{*,n}(md) \right \| ^{2}_{U_{1}})^{\frac{1}{2} } \nonumber
\\&~~~+L_{T}(\left \| \varphi ^{*} \right \|_{h}^{2}+\left \| \varphi ^{\varepsilon} \right \|_{h}^{2} +1)d \nonumber
\\&\le \frac{T}{d} \max_{0\le m\le \left [ (\tau^{*}\wedge T)/d \right ]-1, m\in \mathbb{N}^{+} } (\left \|\mathbb{E}\int_{md}^{(m+1)d} f^{\varepsilon}_{n}(s,u ^{*,n}_{md})-f^{*}_{n}(u ^{*,n}_{md}) \text{d}s\right \|^{2}_{U_{1}})^{\frac{1}{2} }\cdot L_{T}(\left \| \varphi ^{\varepsilon } \nonumber
 \right \|^{2}_{h} +\left \| \varphi ^{*} \right \|^{2}_{h} +1)
 \\&~~~+L_{T}(\left \| \varphi ^{*} \right \|_{h}^{2}+\left \| \varphi ^{\varepsilon} \right \|_{h}^{2} +1)d. \nonumber
	\end{align}
Recalling the definition of the stopping time $\tau^{*}$,  we obtain
\begin{align*}
	|f_{n}(z ,u ^{*,n}_{md})-f(z ,u ^{*,n}_{md})|\le\Re^{\varepsilon}_{f},
	\end{align*}
and
\begin{align*}
	|f^{*}(u ^{*,n}_{md})-f^{*}_{n}(u ^{*,n}_{md})|\le \Re^{\varepsilon}_{f}.
	\end{align*}
Hence, by Lemma 4.2, we obtain the following result:
\begin{align}\label{a18}
	&(\left \|\mathbb{E}\int_{md}^{(m+1)d} f^{\varepsilon}_{n}(s,u ^{*,n}_{md})-f^{*}_{n}(u ^{*,n}_{md}) \text{d}s\right \|^{2}_{U_{1}})^{\frac{1}{2} }\nonumber
\\&=(\left \|\mathbb{E}\int_{md}^{(m+1)d} f_{n}(\frac{s}{\varepsilon } ,u ^{*,n}_{md})-f^{*}_{n}(u ^{*,n}_{md})  \text{d}s\right \|^{2}_{U_{1}})^{\frac{1}{2} }
\\&\le \varepsilon (\left \|\mathbb{E}\int_{\frac{md}{\varepsilon }}^{\frac{(m+1)d}{\varepsilon }} f_{n}(z ,u ^{*,n}_{md})-f^{*}_{n}(u ^{*,n}_{md}) \text{d}z\right \|^{2}_{U_{1}})^{\frac{1}{2} } \nonumber
\\&=\varepsilon (\left \|\mathbb{E}\int_{\frac{md}{\varepsilon }}^{\frac{(m+1)d}{\varepsilon }} f_{n}(z ,u ^{*,n}_{md})-f(z ,u ^{*,n}_{md})+f(z ,u ^{*,n}_{md})-f^{*}(u ^{*,n}_{md})+f^{*}(u ^{*,n}_{md})-f^{*}_{n}(u ^{*,n}_{md}) \text{d}z\right \|^{2}_{U_{1}})^{\frac{1}{2} } \nonumber
\\&\le d[\Phi _{1}(\frac{d}{\varepsilon } )(\left \| \varphi ^{*} \right \| _{h}^{2}+1)+2\Re^{\beta}(n)].\nonumber
	\end{align}
Substituting \eqref{a18} into \eqref{a17} gives
\begin{align}\label{a19}
	\mathfrak{f} ^{2,2}_{2}\le L_{T}(\left \| \varphi ^{\varepsilon }  \right \|^{2}_{h} +\left \| \varphi ^{*} \right \|^{2}_{h} +1)(\Phi _{1}(\frac{d}{\varepsilon } )+d+2\Re^{\beta}(n)).
	\end{align}
Hence substituting \eqref{a9}-\eqref{a18} into \eqref{a8} implies
\begin{align}\label{a19}
	&\mathbb{E}\underset{t\in [0,\tau^{*}\wedge T]}{\sup}\int_{0}^{t}\langle f^{\varepsilon}_{n}(s,u^{\varepsilon,n }_{s})- f^{*}_{n}(u^{*,n }_{s}),u^{\varepsilon,n }(s)- u^{*,n }(s)\rangle_{U_{1}}\text{d}s  \nonumber
\\&\le L_{T,M}[\int_{0}^{ T}\mathbb{E}\underset{z\in [0,\tau^{*}\wedge s]}{\sup}\left \| u^{\varepsilon,n }(z) -u^{*,n }(z) \right \|^{2 }_{U_{1}}\text{d}s+ \left \| \varphi ^{\varepsilon } -\varphi ^{*} \right \| ^{2 }_{h}]
 \\&~~~+L_{T,M}(\left \| \varphi ^{\varepsilon }  \right \|^{2}_{h} +\left \| \varphi ^{*} \right \|^{2}_{h} +1)(\Phi _{1}(\frac{d}{\varepsilon } )+d+d^{\frac{1}{4} }+\Re^{\beta}(n)). \nonumber
	\end{align}

For the diffusion coefficient, we see
\begin{align}\label{a20}
	&\mathbb{E}\int_{0}^{\tau^{*}\wedge T}\left \| g^{\varepsilon}_{n}(s,u^{\varepsilon,n }_{s})- g^{*}_{n}(u^{*,n }_{s})\right \|^{2}_{\mathscr{L}(U_{2},U_{1})}\text{d}s   \nonumber
\\&\le \mathbb{E}\int_{0}^{\tau^{*}\wedge T}\left \| g^{\varepsilon}_{n}(s,u^{\varepsilon,n }_{s})- g^{\varepsilon}_{n}(s,u^{*,n }_{s})\right \|^{2}_{\mathscr{L}(U_{2},U_{1})}\text{d}s
\\&~~~+\mathbb{E}\int_{0}^{\tau^{*}\wedge T}\left \|  g^{\varepsilon}_{n}(s,u^{*,n }_{s})-g^{*}_{n}(u^{*,n }_{s})\right \|^{2}_{\mathscr{L}(U_{2},U_{1})}\text{d}s
\\&:=\mathfrak{g} _{1}+\mathfrak{g}  _{2}.   \nonumber
	\end{align}
By \textbf{(H5)} and \eqref{z1}, we have
\begin{align}\label{a21}
	\mathfrak{g} _{1}&=\mathbb{E}\int_{0}^{\tau^{*}\wedge T}\left \| g^{\varepsilon}_{n}(s,u^{\varepsilon,n }_{s})- g^{\varepsilon}_{n}(s,u^{*,n }_{s})\right \|^{2}_{\mathscr{L}(U_{2},U_{1})}\text{d}s  \nonumber
\\&\le L_{M}[\int_{0}^{ T}\mathbb{E}\underset{z\in [0,\tau^{*}\wedge s]}{\sup}\left \| u^{\varepsilon,n }(s) -u^{*,n }(s) \right \|^{2 }_{U_{1}}\text{d}s+\left \| \varphi ^{\varepsilon } -\varphi ^{*} \right \| ^{2 }_{h}],
	\end{align}
and
\begin{align}\label{a22}
	\mathfrak{g} _{2}&=\mathbb{E}\int_{0}^{\tau^{*}\wedge T}\left \|  g^{\varepsilon}_{n}(s,u^{*,n }_{s})-g^{*}_{n}(u^{*,n }_{s})\right \|^{2}_{\mathscr{L}(U_{2},U_{1})}\text{d}s  \nonumber
\\&\le\mathbb{E}\int_{0}^{\tau^{*}\wedge T}\left \|  g^{\varepsilon}_{n}(s,u^{*,n }_{s})-g^{\varepsilon}_{n}(s,\hat{u}^{*,n }_{s})\right \|^{2}_{\mathscr{L}(U_{2},U_{1})}\text{d}s\nonumber
\\&~~~+\mathbb{E}\int_{0}^{\tau^{*}\wedge T}\left \|  g^{\varepsilon}_{n}(s,\hat{u}^{*,n }_{s})-g^{*}_{n}(\hat{u}^{*,n }_{s})\right \|^{2}_{\mathscr{L}(U_{2},U_{1})}\text{d}s
\\&~~~+\mathbb{E}\int_{0}^{\tau^{*}\wedge T}\left \|  g^{*}_{n}(\hat{u}^{*,n }_{s})-g^{*}_{n}(u^{*,n }_{s})\right \|^{2}_{\mathscr{L}(U_{2},U_{1})}\text{d}s \nonumber
\\&:=\mathfrak{g} ^{2,1}_{2}+\mathfrak{g} ^{2,2}_{2}+\mathfrak{g} ^{2,3}_{2}. \nonumber
	\end{align}
Then by \textbf{(H4)}, \eqref{p7}, H\"older's inequality, Jensen's inequality and Lemma 4.2, we have
\begin{align}\label{a23}
	\mathfrak{g}^{2,1} _{2}&=\mathbb{E}\int_{0}^{\tau^{*}\wedge T}\left \|  g^{\varepsilon}_{n}(s,u^{*,n }_{s})-g^{\varepsilon}_{n}(s,\hat{u}^{*,n }_{s})\right \|^{2}_{\mathscr{L}(U_{2},U_{1})}\text{d}s \nonumber
\\&\le L_{T}\mathbb{E}\int_{0}^{\tau^{*}\wedge T}\left \|  u^{*,n }_{s}-\hat{u}^{* ,n}_{s}\right \|^{2}_{h}\text{d}s
\\&\le L_{T}(1+\left \| \varphi^{*} \right \|_{h}^{2} )d^{\frac{1 }{2} }
. \nonumber
	\end{align}
Similarly, for $\mathfrak{g}^{2,3} _{2}$,
\begin{align}\label{a25}
	\mathfrak{g}^{2,3} _{2}&=\mathbb{E}\int_{0}^{\tau^{*}\wedge T}\left \|  g^{*}_{n}(\hat{u}^{* ,n}_{s})-g^{*}_{n}(u^{*,n }_{s})\right \|^{2}_{\mathscr{L}(U_{2},U_{1})}\text{d}s
\\&\le L_{T}(1+\left \| \varphi^{*} \right \|_{h}^{2} )d^{\frac{1 }{2} }
. \nonumber
	\end{align}
Then, we will use the time discretization technique to deal with $\mathfrak{g}^{2,2} _{2}$:
\begin{align}\label{a26}
	\mathfrak{g}^{2,2} _{2}&=\mathbb{E}\int_{0}^{\tau^{*}\wedge T}\left \|  g^{\varepsilon}_{n}(s,\hat{u}^{*,n }_{s})-g^{*}_{n}(\hat{u}^{*,n }_{s})\right \|^{2}_{\mathscr{L}(U_{2},U_{1})}\text{d}s   \nonumber
\\&=\mathbb{E}\sum_{n=0}^{\left [ \frac{\tau^{*}\wedge T}{d}  \right ]-1 } \int_{nd}^{(n+1)d}\left \|  g^{\varepsilon}_{n}(s,u^{*,n }_{nd})-g^{*}_{n}(u^{*,n }_{nd})\right \|^{2}_{\mathscr{L}(U_{2},U_{1})}\text{d}s  \nonumber
\\&~~~+\mathbb{E}\int_{\left [ \frac{\tau^{*}\wedge T}{d}  \right ]d}^{\tau^{*}\wedge T}\left \|  g^{\varepsilon}_{n}(s,u^{*,n }_{\left [ \tau^{*}\wedge T/d  \right ]d})-g^{*}_{n}(u^{*,n }_{\left [ \tau^{*}\wedge T/d  \right ]d})\right \|^{2}_{\mathscr{L}(U_{2},U_{1})}\text{d}s   \nonumber
\\& \le \varepsilon\sum_{n=0}^{\left [ \frac{\tau^{*}\wedge T}{d}  \right ]-1 } \mathbb{E}\int_{\frac{nd}{\varepsilon }}^{\frac{(n+1)d}{\varepsilon }}\left \|  g(z,u^{*,n }_{nd})-g^{*}_{n}(u^{*,n }_{nd})\right \|^{2}_{\mathscr{L}(U_{2},U_{1})}\text{d}z
\\&~~~+L_{M}(1+\mathbb{E}\left \| u^{*,n }_{\left [ \tau^{*}\wedge T/d  \right ]d}\right \|^{2}_{h})d   \nonumber
\\&\le \sum_{n=0}^{\left [ \frac{\tau^{*}\wedge T}{d}  \right ]-1 }d[\Phi _{2}(\frac{d}{\varepsilon } )(\left \| \varphi ^{*} \right \| _{h}^{2}+1)+2\Re^{\gamma}(n)]+L_{T}(\left \| \varphi ^{*} \right \| _{h}^{2}+1)d  \nonumber
\\&\le L_{M,T}(\left \| \varphi ^{*} \right \| _{h}^{2}+1)[d+\Phi _{2}(\frac{d}{\varepsilon } )+\Re^{\gamma}(n)]
. \nonumber
	\end{align}
In conclusion, by \eqref{a21}-\eqref{a26}, we obtain
\begin{align}\label{a27}
	&\mathbb{E}\int_{0}^{\tau^{*}\wedge T}\left \| g^{\varepsilon}_{n}(s,u^{\varepsilon,n }_{s})- g^{*}_{n}(u^{*,n }_{s})\right \|^{2}_{\mathscr{L}(U_{2},U_{1})}\text{d}s   \nonumber
\\&\le L_{T,M}[\int_{0}^{ T}\mathbb{E}\underset{z\in [0,\tau^{*}\wedge s]}{\sup}\left \| u^{\varepsilon,n }(z) -u^{*,n }(z) \right \|^{2 }_{U_{1}}\text{d}s+ \left \| \varphi ^{\varepsilon } -\varphi ^{*} \right \| ^{2 }_{h}]
 \\&~~~+L_{T,M}(\left \| \varphi ^{\varepsilon }  \right \|^{2}_{h} +\left \| \varphi ^{*} \right \|^{2}_{h} +1)[\Phi _{2}(\frac{d}{\varepsilon } )+d+d^{\frac{1 }{2} }+\Re^{\gamma}(n)]. \nonumber
	\end{align}
Let $d=\varepsilon ^{\frac{1}{2} }$, then substituting \eqref{a17}-\eqref{a27} into \eqref{a7} implies
\begin{align*}
		&\mathbb{E}\underset{t\in [0, T]}{\sup}\left \| u^{\varepsilon,n }(\tau^{*}\wedge t;\varphi ^{\varepsilon }) -u^{*,n }(\tau^{*}\wedge t;\varphi ^{*}) \right \|^{2}_{U_{1}}  \nonumber
\\&\le L_{T,M}[\left \| \varphi ^{\varepsilon } -\varphi ^{*} \right \|^{2}_{h} +\int_{0}^{ T}\mathbb{E}\underset{z\in [0, s]}{\sup}\left \| u^{\varepsilon,n }(\tau^{*}\wedge z)- u^{*,n }(\tau^{*}\wedge z) \right \| _{U_{1}}^{2 }\text{d}s]  \nonumber
\\&~~~+L_{T,M}(\left \| \varphi ^{\varepsilon }  \right \|^{2}_{h} +\left \| \varphi ^{*} \right \|^{2}_{h} +1)[\Phi _{1}(\varepsilon ^{-\frac{1}{2} } )+\Phi _{2}(\varepsilon ^{-\frac{1}{2} } )+\varepsilon^{\frac{1 }{2} }+\varepsilon^{\frac{1 }{8} }+\varepsilon^{\frac{1}{4} }+\Re^{\beta}(n)+\Re^{\gamma}(n)]. \nonumber
	\end{align*}
It follows from Gronwall's lemma that
\begin{align}\label{a27aa}
	&\mathbb{E}\underset{t\in [0, \tau^{*}\wedge T]}{\sup}\left \| u^{\varepsilon,n }( t;\varphi ^{\varepsilon }) -u^{*,n }(t;\varphi ^{*}) \right \|^{2}_{U_{1}}\nonumber
\\&\le L_{T,M}[\left \| \varphi ^{\varepsilon } -\varphi ^{*} \right \|^{2}_{h} +(\left \| \varphi ^{\varepsilon }  \right \|^{2}_{h} +\left \| \varphi ^{*} \right \|^{2}_{h} +1)(\Phi _{1}(\varepsilon ^{-\frac{1}{2} } )+\Phi _{2}(\varepsilon ^{-\frac{1}{2} } )
\\&~~~+\varepsilon^{\frac{1 }{8}  }+\Re^{\beta}(n)+\Re^{\gamma}(n))]\exp\{L_{T,M}\}.\nonumber
	\end{align}

Similar to \eqref{eq:condition}-\eqref{i15i}, we obtain that  for any $\epsilon>0$, there exists a compact subset $\mathcal{C}\subseteq C^{h}_{U_{1}}$ such that
\begin{align*}
\mathbb{P}\left( u_{t}^{\varepsilon}\in \mathcal{C}~ \text{and} ~ u_{t}^{*}\in  \mathcal{C}\right) \geq 1 - \epsilon, \quad \forall t \in [0,T],
\end{align*}
 which implies  the processes $u_{t}^{\varepsilon}$ and $u_{t}^{*}$ remain confined within $\mathcal{C}$ for all $t\in[0, T]$ with probability at least  $1-\epsilon$. Consequently, we can choose a compact set $\mathcal{C}$ such that
\begin{align}\label{i999i}
\mathbb{P}(\tau^{*} \leq T)<\epsilon.
\end{align}
Hence, when $ T\leq\tau^{*} $, provided that $n>N$, it follows that
\begin{align*}
&\mathbb{E}\underset{t\in [0,  T]}{\sup}\left \| u^{\varepsilon,n }( t;\varphi ^{\varepsilon }) -u^{*,n }(t;\varphi ^{*}) \right \|^{2}_{U_{1}}
\\&\le L_{T,M}\sqrt{[\left \| \varphi ^{\varepsilon } -\varphi ^{*} \right \|^{2}_{h} +(\left \| \varphi ^{\varepsilon }  \right \|^{2}_{h} +\left \| \varphi ^{*} \right \|^{2}_{h} +1)(\Phi _{1}(\varepsilon ^{-\frac{1}{2} } )+\Phi _{2}(\varepsilon ^{-\frac{1}{2} } )
+\varepsilon^{\frac{1 }{8}  }+\varepsilon^{\beta}+\varepsilon^{\gamma})]},\nonumber
\end{align*}
which implies that for any   real-valued, bounded and continuous function $\mathcal{L}$ defined on the space $C([0, T], U_{1})$,
\begin{align}\label{c999c}
\lim_{\varepsilon  \to 0}|\mathbb{E}\mathcal{L}(u^{\varepsilon,n}(t)| _{[0,T]}) - \mathbb{E}\mathcal{L}(u^{*,n}(t)|_{[0,T]})|=0.
\end{align}
Finally, we follow the proof of weak uniqueness given in Theorem 3.2 and consequently obtain
\begin{align*}
\mathbb{E}\mathcal{L}(u^{\varepsilon,n}(t)| _{[0,T]}) - \mathbb{E}\mathcal{L}(u^{\varepsilon}(t)|_{[0,T]})\to 0,\quad n\to \infty,
\end{align*}
and
\begin{align*}
\mathbb{E}\mathcal{L}(u^{*,n}(t)| _{[0,T]}) - \mathbb{E}\mathcal{L}(u^{*}(t)|_{[0,T]})\to 0,\quad n\to \infty,
\end{align*}
which implies that for any $\varepsilon>0$, there exists an integer $N$ such that for all $n>N$,
\begin{align}\label{a999a}
|\mathbb{E}\mathcal{L}(u^{\varepsilon,n}(t)| _{[0,T]}) - \mathbb{E}\mathcal{L}(u^{\varepsilon}(t)|_{[0,T]})|<\varepsilon,
\end{align}
and
\begin{align}\label{b999b}
|\mathbb{E}\mathcal{L}(u^{*,n}(t)| _{[0,T]}) - \mathbb{E}\mathcal{L}(u^{*}(t)|_{[0,T]})|<\varepsilon,
\end{align}
then, applying \eqref{c999c}-\eqref{b999b}, we obtain that when $ T\leq\tau^{*} $,
\begin{align*}
\lim_{\varepsilon  \to 0}|\mathbb{E}\mathcal{L}(u^{\varepsilon}(t)| _{[0,T]}) - \mathbb{E}\mathcal{L}(u^{*}(t)|_{[0,T]})|=0.
\end{align*}
Combining \eqref{i999i} and letting $\epsilon\to 0$, we obtain
\begin{align}\label{c1}
		\lim_{\varepsilon  \to 0}|\mathbb{E}\mathcal{L}(u^{\varepsilon }(t;\varphi ^{\varepsilon })| _{[0,T]}) - \mathbb{E}\mathcal{L}(u^{* }(t;\varphi ^{*})|_{[0,T]})|= 0.
	\end{align}
This completes the proof.\quad $\Box $
\\~
\\\textbf{Remark 4.4.} Similarly, we emphasize that the methods and conclusions in the paper can also be applied to the corresponding models with finite delay and without delay, with minor modifications to the conditions.
~\\

Consider the following SFPDEs with finite delay
\begin{align}\label{j7}
		\begin{cases}
 \text{d}u^{\varepsilon }(t)=(A(u^{\varepsilon }(t))+f(\frac{t}{\varepsilon } ,u^{\varepsilon }_{t})\text{d}t+g(\frac{t}{\varepsilon } ,u^{\varepsilon }_{t})\text{d}W(t), \\
u_{0} =\varphi^{\varepsilon }\in \mathcal{H} ,
\end{cases}
	\end{align}
where $\varepsilon \in (0,1]$.
~\\
\\\textbf{(H6')} There exist functions $\Phi _{1} $, $\Phi _{2} $ and $f^{*}\in C(\mathcal{H},U_{1})$, $g^{*}\in C(\mathcal{H},\mathscr{L}(U_{2},U_{1}) )$ such that for any $t\ge0$, $r>0$ and $\varphi \in \mathcal{H}$,
\begin{align*}
		\frac{1}{r} \left \| \int_{t}^{t+r}[f(s,\varphi ) -f^{*}(\varphi )] \text{d}s\right \| _{U_{1}}\le \Phi _{1}(r)(\left \| \varphi  \right \|_{h}+M ),
	\end{align*}
\begin{align*}
		\frac{1}{r}\int_{t}^{t+r} \left \| g(s,\varphi ) -g^{*}(\varphi )\right \|^{2} _{\mathscr{L}(U_{2},U_{1}) } \text{d}s\le \Phi_{2} (r)(\left \| \varphi  \right \|^{2}_{h}+M ).
	\end{align*}
For the following averaged equation
\begin{align}\label{j8}
		\begin{cases}
 \text{d}u(t)=(A(u(t))+f^{*}_{n}(u_{t})\text{d}t+g^{*}_{n}(u_{t})\text{d}W(t), \\
u_{0} =\varphi^{* }\in \mathcal{H},
\end{cases}
	\end{align}
we have the following theorem:
~\\
\\\textbf{Theorem 4.5.} \emph{Consider \eqref{j7} and \eqref{j8}. Suppose that the assumptions \textbf{(H1)}$-$\textbf{(H5)} and \textbf{(H6')} hold. For any initial values $\varphi ^{\varepsilon }, \varphi ^{*}\in \mathcal{H}$ and $T > 0$, assume further that $$\lim_{\varepsilon  \to 0} \mathbb{E} \left \|\varphi ^{\varepsilon } -\varphi ^{*} \right \|^{2}_{h}=0,$$ then we have
\begin{align}\label{c1}
		\lim_{\varepsilon  \to 0}|\mathbb{E}\mathcal{L}(u^{\varepsilon }(t;\varphi ^{\varepsilon })| _{[0,T]}) - \mathbb{E}\mathcal{L}(u^{* }(t;\varphi ^{*})|_{[0,T]})|= 0,
	\end{align}
where  $u^{\varepsilon}(s;\varphi ^{\varepsilon })$ is the weak solution of \eqref{j7} and
$u^{* }(s;\varphi ^{*})$ is the weak solution of \eqref{j8}.}
~\\

For the following SPDEs:
\begin{align}\label{j9}
		\begin{cases}
 \text{d}u^{\varepsilon }(t)=(A(u^{\varepsilon }(t))+f(\frac{t}{\varepsilon } ,u^{\varepsilon }(t))\text{d}t+g(\frac{t}{\varepsilon } ,u^{\varepsilon }(t))\text{d}W(t), \\
u(0) =u^{\varepsilon }_{0}\in U_{1} ,
\end{cases}
	\end{align}
let $\varepsilon \in (0,1]$.
~\\
\\\textbf{(H7)} There exist functions $\Phi _{1} $, $\Phi _{2} $ and $f^{*}\in C(U_{1},U_{1})$, $g^{*}\in C(U_{1},\mathscr{L}(U_{2},U_{1}) )$ such that for any $t\ge0$, $r>0$ and $u \in U_{1}$,
\begin{align*}
		\frac{1}{r} \left \| \int_{t}^{t+r}[f(s,u ) -f^{*}(u )] \text{d}s\right \| _{U_{1}}\le \Phi _{1}(r)(\left \| u  \right \|_{U_{1}}+M ),
	\end{align*}
\begin{align*}
		\frac{1}{r}\int_{t}^{t+r} \left \| g(s,u ) -g^{*}(u )\right \|^{2} _{\mathscr{L}(U_{2},U_{1}) } \text{d}s\le \Phi_{2} (r)(\left \| u \right \|^{2}_{U_{1}}+M ).
	\end{align*}
For the following averaged equation
\begin{align}\label{j10}
		\begin{cases}
 \text{d}u(t)=(A(u(t))+f^{*}(u(t))\text{d}t+g^{*}(u(t))\text{d}W(t), \\
u(0) =u^{*}_{0}\in U_{1},
\end{cases}
	\end{align}
we have the following theorem:
~\\
\\\textbf{Theorem 4.6.} \emph{Consider \eqref{j7} and \eqref{j8}. Suppose that the assumptions \textbf{(H1)}, \textbf{(H2')}, \textbf{(H3)}, \textbf{(H4')} and \textbf{(H7)} hold.  For any initial values $u ^{\varepsilon }_{0}, u ^{*}_{0}\in U_{1}$ and $T > 0$, assume further that $\lim_{\varepsilon  \to 0} \mathbb{E} \left \|u ^{\varepsilon }_{0}-u ^{*}_{0} \right \|^{2}_{U_{1}}=0$. Then we have
\begin{align*}
		\lim_{\varepsilon  \to 0}|\mathbb{E}\mathcal{L}(u^{\varepsilon }(t;u ^{\varepsilon }_{0})| _{[0,T]}) - \mathbb{E}\mathcal{L}(u^{* }(t;u ^{*}_{0})|_{[0,T]})|= 0,
	\end{align*}}
where  $u^{\varepsilon}(t;u ^{\varepsilon }_{0})$ is the weak solution of \eqref{j9} and
$u^{*}(t; u ^{*}_{0})$ is the weak solution of \eqref{j10}.

The proofs of Theorems 4.5 and 4.6 follow a procedure analogous to the previously discussed proofs of Theorem 4.3; therefore, we omit these particular details.
\section{\textup{Applications}}
This section aims to substantiate the validity of our principal findings by applying them to stochastic generalized porous media equations and stochastic reaction diffusion equations.  Consider $D\subset \mathbb{R}^{n}(n \in \mathbb{N}) $ as an open bounded subset, with $ - \Delta$ governed by Dirichlet boundary conditions.

Consider the stochastic generalized porous media equation:
\begin{equation}\label{g31}
\begin{cases}
\begin{aligned}
 \text{d}u&=[\Delta(\left | u \right | ^{q-2}+u) +\xi_{1} (\frac{t}{\varepsilon } )f(u)]\text{d}t+\xi_{2} (\frac{t}{\varepsilon } )g(u)\text{d}W(t),
 \end{aligned}\\
 u(0)=u_{0}^{\varepsilon}\in W^{1,2}_{0}(D ),
\end{cases}.
\end{equation}
where  $0<\varepsilon \ll 1$, $W$ is a  standard real-valued Wiener process and $q>2$. We require both functions $\xi_{1}$ and $\xi_{2}$ to be positively bounded, meaning there exists a constant $M$ such that $\left | \xi _{i}(t) \right | \le M$ for any $t\in \mathbb{R} ^{+}$ and $i=1,2$.

Assume that $f$ and $g$ satisfy (\textbf{H2'}) and (\textbf{H4'}). Then, by Theorem 3.7, we obtain the existence of a  unique weak solution  $u^{\varepsilon}(t;u_{0}^{\varepsilon})$ for system \eqref{g31}, i.e., $u^{\varepsilon}(t;u_{0}^{\varepsilon})$ is unique in  law. In addition,  we can also obtain the following theorem:
~\\
\\\textbf{Theorem 5.1.} \emph{Let $u^{*}$ be the  stationary solution of the following averaged equation:
\begin{eqnarray*}
\begin{cases}
\begin{aligned}
 \text{d}u&=[\Delta(\left | u \right | ^{q-2}+u) +\xi_{1} ^{*}f(u)]\text{d}t+\xi_{2}^{*}g(u)\text{d}W(t),
 \end{aligned}\\
 u(0)=u_{0}^{*}\in W^{1,2}_{0}(D ),
\end{cases}.
\end{eqnarray*}
where $\xi_{i}^{*}=\lim_{T \to \infty} \frac{1}{T} \int_{t}^{t+T} \xi_{i}(s)\text{d}s$ for any $t\in \mathbb{R} ^{+}$ and $i=1,2$. Then, assume further that $\lim_{\varepsilon  \to 0}  \left \|u_{0} ^{\varepsilon} -u_{0}^{*} \right \|^{2}_{W^{1,2}_{0}}=0$, we have
\begin{align*}
		\lim_{\varepsilon  \to 0}|\mathbb{E}\mathcal{L}(u^{\varepsilon }(t;u ^{\varepsilon }_{0})| _{[0,T]}) - \mathbb{E}\mathcal{L}(u^{* }(t;u ^{*}_{0})|_{[0,T]})|= 0,
	\end{align*}
for any $T>0$.
~\\
\\\textbf{proof}} Let $B=L^{q}(D)$, $U_{1}=W^{1,2}_{0}(D)$ and
\begin{eqnarray*}
\begin{split}
		_{V ^{\ast}} \langle A(u), v\rangle _{V}:=-\int _{D }u(x)\left | u(x) \right |^{q-2}v(x)\text{d}x-a\int _{D }u(x)v(x)\text{d}x,
\end{split}
	\end{eqnarray*}
for $u,v\in B$, which implies that $B\subset  U_{1}= U_{1}^{*}\subset B^{*}$. Then, according to Theorem 6.3 in Reference \cite{ref2}, the operator $A(u)=\Delta(\left | u \right | ^{q-2}+u)$ satisfies conditions (\textbf{H1}), (\textbf{H2'}) and (\textbf{H3}). In fact, the operator $A$ satisfies the monotonicity condition. Therefore, according to Theorem 4.6, we can derive the the first Bogolyubov theorem.

For the specific form of the diffusion and drift coefficients, for example, let
\begin{eqnarray*}
\begin{split}
		f(u)= \sin \sqrt{ \left | u \right |  } ,\quad g(u)=\cos (u^{\frac{2}{3}}),
\end{split}
	\end{eqnarray*}
for any $u\in U_{1}$. Then, it is easy to check that $f$ and $g$ satisfy the Assumption (\textbf{H4'}) with H$\ddot{\text{o}}$lder exponent $\beta=\frac{1}{2} $ and $\gamma =\frac{2}{3}$.
~\\

For the stochastic reaction diffusion equations with infinite delay:
\begin{eqnarray*}
\begin{cases}
\begin{aligned}
 \text{d}u^{\varepsilon}(t)&=[\Delta u^{\varepsilon} -u^{\varepsilon}\left | u^{\varepsilon} \right | ^{q-2}+\xi_{1} (\frac{t}{\varepsilon } )f(u^{\varepsilon}_{t})]\text{d}t+g\text{d}W(t),
 \end{aligned}\\
 u_{0}=\varphi^{\varepsilon}\in C^{h}_{L^{2}(D )},
\end{cases}.
\end{eqnarray*}
let $W$ is a one-dimensional two-sided  cylindrical $Q$-Wiener process with $Q = I$ on $B=H^{1,2}_{0}(D)\cap L^{q}(D)$ and $g \in \mathscr{L}(B,L^{2}(D))$. Assume that $f$ satisfies (\textbf{H2}),  (\textbf{H4}) and $g$ satisfies (\textbf{H5}), then we have the following theorem:
~\\
\\\textbf{Theorem 5.2.} \emph{Let $u^{*}$ be the  stationary solution of the following averaged equation:
\begin{eqnarray*}
\begin{cases}
\begin{aligned}
 \text{d}u(t)&=[\Delta u -u\left | u \right | ^{q-2} +\xi_{1} ^{*}f(u_{t})]\text{d}t+g\text{d}W(t),
 \end{aligned}\\
 u_{0}=\varphi^{*}\in C^{h}_{L^{2}(D )}.
\end{cases}.
\end{eqnarray*}
 Then, assume further that $\lim_{\varepsilon  \to 0} \mathbb{E} \left \|\varphi ^{\varepsilon } -\varphi ^{*} \right \|^{2}_{h}=0$, we have
\begin{align*}
		\lim_{\varepsilon  \to 0}|\mathbb{E}\mathcal{L}(u^{\varepsilon }(t;\varphi ^{\varepsilon })| _{[0,T]}) - \mathbb{E}\mathcal{L}(u^{* }(t;\varphi ^{*})|_{[0,T]})|= 0,
	\end{align*}
for any $T>0$.
~\\
\\\textbf{proof}} Let $B=H^{1,2}_{0}(D)\cap L^{q}(D)$, $U_{1}=L^{2}(D)$ and according to Theorem 6.1 in \cite{ref2}, the operator $A(u)=\Delta u -u\left | u \right | ^{q-2}$ satisfies conditions (\textbf{H1})-(\textbf{H3}). Therefore, according to Theorem 4.3, we can derive the aforementioned conclusion.

For the specific form of $f$, let
\begin{eqnarray*}
\begin{split}
		f(\varphi)=\cos \sqrt{\left | \varphi (0) \right | }  +\int_{-\infty    }^{0} \sqrt{\left | \varphi (\theta ) \right | } \mu (\text{d}\theta ) ,
\end{split}
	\end{eqnarray*}
for any $\varphi\in C^{h}_{L^{2}(D )}$, where $\mu (\text{d}\theta)=2he^{2h\theta }\text{d}\theta $. Then, it is easy to check that $f$ satisfies  the Assumption (\textbf{H2}) and  (\textbf{H4}) with H$\ddot{\text{o}}$lder exponent $\beta =\frac{1}{2} $.


\section*{Acknowledgments}
 The first author (S. Lu)  supported by Graduate Innovation Fund of Jilin University. The second author (X. Yang) was supported by  National Natural Science Foundation of China (12071175, 12371191). The third author (Y. Li) was supported by  National Natural Science Foundation of China (12071175 and 12471183).
\section*{Data availability}
No data was used for the research described in the article.

\section*{References}
\bibliographystyle{plain}
\bibliography{ref}

\end{document}